\documentclass{amsart}

\title{Vanishing Integrals for Hall--Littlewood polynomials}
\author{Vidya Venkateswaran}
\address{Department of Mathematics, California Institute of Technology, Pasadena, CA 91125}
\email{\href{mailto:vidyav@caltech.edu}{vidyav@caltech.edu}}
\usepackage{hyperref}
\usepackage[svgnames]{xcolor}
\hypersetup{colorlinks,breaklinks,
            linkcolor=DarkBlue,urlcolor=DarkBlue,
            anchorcolor=DarkBlue,citecolor=DarkBlue}
\usepackage{url}
\usepackage{euscript}
\usepackage[verbose,letterpaper,tmargin=1in,bmargin=1in,lmargin=1in,rmargin=1in]{geometry}
\usepackage{amssymb}
\usepackage{amsmath}

\usepackage{verbatim}

\numberwithin{equation}{section}

\newtheorem{theorem}{Theorem}[section]
\newtheorem{claim}{Claim}[theorem]

\newtheorem{lemma}[theorem]{Lemma}
\newtheorem{proposition}[theorem]{Proposition}
\newtheorem{corollary}[theorem]{Corollary}

\theoremstyle{definition}

\theoremstyle{remark}

\DeclareMathOperator{\PV}{P.V.}

\newcommand{\qbinom}[2]{\left[#1 \atop #2\right]}

\begin{document}

\begin{abstract}
It is well known that if one integrates a Schur function indexed by a partition $\lambda$ over the symplectic (resp. orthogonal) group, the integral vanishes unless all parts of $\lambda$ have even multiplicity (resp. all parts of $\lambda$ are even).  In a recent paper of Rains and Vazirani, Macdonald polynomial generalizations of these identities and several others were developed and proved using Hecke algebra techniques.  However at $q=0$ (the Hall--Littlewood level), these approaches do not work, although one can obtain the results by taking the appropriate limit.  In this paper, we develop a direct approach for dealing with this special case.  This technique allows us to prove some identities that were not amenable to the Hecke algebra approach.  Moreover, we are able to generalize some of the identities by introducing extra parameters.  This leads us to a finite-dimensional analog of a recent result of Warnaar, which uses the Rogers--Szeg\H{o} polynomials to unify some existing summation type formulas for Hall--Littlewood functions.  
\end{abstract}

\maketitle

\section{Introduction}
Two classical identities in the representation theory of real Lie groups are:
\begin{theorem}
For any integer $n \geq 0$ and partition $\lambda$ with at most $n$ parts, we have 
\begin{align*}
\int_{O \in O(n)} s_{\lambda}(O) dO
= \begin{cases} 1, &\text{if all parts of $\lambda$ are even } \\
0, & \text{otherwise}
\end{cases}
\end{align*}
(where the integral is with respect to Haar measure on the orthogonal group).  Similarly, for $n$ even, we have
\begin{align*}
\int_{S \in Sp(n)} s_{\lambda}(S) dS
=\begin{cases} 1, &\text{if all parts of $\lambda$ have even multiplicity} \\
0, & \text{otherwise}
\end{cases}
\end{align*}
(where the integral is with respect to Haar measure on the symplectic group).
\end{theorem}
Here $s_{\lambda}$ is the Schur function in $n$ variables indexed by the partition $\lambda$.  Schur functions have an intimate connection to representation theory: they give the character of an irreducible representation of the unitary group, $U(n)$.  In particular, the character's value on a matrix is given by evaluating the Schur function at the matrix's eigenvalues.  Thus, the above identities encode the following facts: in the expansion of $s_{\lambda}$ into irreducible characters of $O(n)$ (resp. $Sp(n)$), the coefficient of the trivial character is zero unless all parts of $\lambda$ are even (resp. all parts of $\lambda$ have even multiplicity).  These identities can be proved using the Gelfand pairs $(G,K) = (GL_{n}(\mathbb{R}), O(n))$ and $(GL_{n}(\mathbb{H}), U(n, \mathbb{H}))$ and the decomposition of the induced trivial representation into irreducible representations of $G$, see \cite{Mac}.  For example, the orthogonal group identity follows from the structure result
\begin{align*}
e_{K} P(G) = P(K \backslash G) \cong \displaystyle \bigoplus_{l(\lambda) \leq n} F_{2\lambda}(V)
\end{align*}
(in the notation of \cite{Mac}) and the fact that the Schur function gives the character of a polynomial representation of $GL_{n}(\mathbb{R})$.

Note that using the eigenvalue densities for the orthogonal and symplectic groups, we may rephrase the above identities in terms of random matrix averages.  For example, the symplectic integral above can be rephrased as
\begin{align*}
\frac{1}{2^{n}n!}\int_{T} s_{\lambda}(z_{1}, z_{1}^{-1}, z_{2}, z_{2}^{-1}, \dots, z_{n}, z_{n}^{-1}) \prod_{1 \leq i \leq n} |z_{i} - z_{i}^{-1}|^{2} \prod_{1 \leq i<j \leq n} |z_{i} + z_{i}^{-1} - z_{j} - z_{j}^{-1}|^{2} dT,
\end{align*}
where 
\begin{align*}
T &= \{ (z_{1}, \dots, z_{n}) : |z_{1}| = \dots = |z_{n}| = 1 \} \\
dT &= \prod_{j} \frac{dz_{j}}{2 \pi \sqrt{-1} z_{j}}
\end{align*}
are the $n$-torus and Haar measure, respectively.  Such identities, and their generalizations, have consequences outside symmetric function theory.  For example, in their work dealing with symmetrized generalizations of the Hammersely process \cite{FR}, Forrester and Rains developed an $\alpha$-generalization of the above orthogonal group identity.  

A natural question, then, is whether such identities admit a $q,t$ generalization to the level of Macdonald polynomials.  In \cite{R}, a number of such identities were conjectured: that is, a suitable choice of density was suggested so that integrating Macdonald polynomials against it should vanish  unless the partition is of the appropriate form, and such that when $q=t$, these identities become the known ones for Schur functions.  In \cite{RV}, Rains and Vazirani developed Hecke algebra techniques which enabled them to prove many of these results.  In fact, only Conjectures 3 and 5 of \cite{R} remain open.  

An interesting subfamily of the Macdonald polynomials are the Hall--Littlewood polynomials which are obtained at $q=0$, see Chapter 3 of \cite{Mac}.  Unfortunately, none of the above proofs work at this level: they involve $q$-shift operators, which do not behave well at $q=0$.  In this paper, we develop a direct method for dealing with these cases.  This method allows us to explicitly obtain the nonzero values as well as generalizations involving extra parameters.  We also prove Conjectures 3 and 5 from \cite{R} for Hall--Littlewood polynomials.  

There are several nice consequences.  The first involves a (recent) identity discovered by Warnaar for Hall--Littlewood polynomials \cite{W}.  He uses the Rogers--Szeg\H{o} polynomials to unify the Littlewood/Macdonald identities for Hall--Littlewood functions.  We find a two-parameter integral identity and, using a method of Rains, we show that in the limit $n \rightarrow \infty$ it becomes Warnaar's identity.  Thus, our identity may be viewed as a finite-dimensional analog of Warnaar's summation result.  Another unexpected feature of the direct method we employ is an underlying Pfaffian structure in the orthogonal cases.  It turns out that Pfaffians of suitable matrices nicely enumerate the nonzero values of these integrals.  While Pfaffians are very common in Schur function identities (Schur functions are ratios of determinants), to our knowledge this is the first time they are appearing in the Hall--Littlewood context.  Finally, the identities below involve Hall--Littlewood polynomials with a parameter $t$, but in many instances the evaluation of the integral produces a polynomial in $t^{2}$ or $\sqrt{t}$ (see for example, the symplectic and Kawanaka integrals, Corollary 6.3 and 6.4 respectively).  Thus these identities may be viewed as quadratic transformations of Hall--Littlewood polynomials.

The outline of the paper is as follows.  In the second section, we introduce some basic notation and review Hall--Littlewood polynomials.  In the third section, we prove Hall--Littlewood orthogonality to illustrate our method of proof in a basic case.  In the next section, we use Pfaffians and some technical arguments to prove $\alpha$ generalizations of the orthogonal integrals.  In section 5, we use a Pieri rule to add one more parameter $\beta$ to these identities.  In section 6, we discuss special cases of the $\alpha, \beta$ identity.  In section 7, we prove that in the limit $n \rightarrow \infty$ of the $\alpha, \beta$ identity, we  recover Warnaar's identity.  Finally, in the last section, we prove some remaining vanishing results from \cite{RV} and \cite{R}.

We mention some related work in progress.  As was discussed above, many of the integral identities for $t=0$ (the Schur case) follow from the theory of symmetric spaces, and thus have a representation theoretic significance.  In \cite{MacP}, Macdonald shows that Hall--Littlewood polynomials (and their analogs for other classical root systems) arise as zonal spherical functions on $p$-adic reductive groups.  Given this construction, it is natural to wonder whether our identities can be interpreted as $p$-adic analogs of the Schur cases.  In a follow-up project, we will show that this is indeed the case: we give another proof via integrals over $p$-adic groups.

Finally, many of the integral identities of \cite{RV} involve Koornwinder polynomials, a $6$-parameter $BC_{n}$-symmetric family of Laurent polynomials that contain the Macdonald polynomials as suitable limits of the parameters.  Just as in the Macdonald polynomial case, standard constructions via difference operators do not allow one to control the $q=0$ polynomials.  The first step in obtaining an analog of the Hall--Littlewood polynomials is to produce a $q=0$ closed form.  Such a formula is not known; in further work we will use orthogonality of Koornwinder polynomials to provide an explicit closed form \cite{VV}.  We then use this result to prove the $q=0$ cases of the remaining identities in \cite{RV}.

\medskip
\noindent\textbf{Acknowledgements.}  The author would like to thank her advisor, E. Rains, for suggesting the problems in this paper and for all the guidance and support he gave her during this project.  She would also like to thank A. Borodin, P. Diaconis, P. Forrester, M. Vazirani, and O. Warnaar for many useful discussions and comments.  

\section{Background and Notation}

We will briefly review Hall--Littlewood polynomials; we follow \cite{Mac}.  We also set up the required notation.

Let $\lambda = (\lambda_{1}, \dots, \lambda_{n})$ be a partition, in which some of the $\lambda_{i}$ may be zero.  In particular, note that $\lambda_{1} \geq \lambda_{2} \geq\cdots \geq \lambda_{n} \geq 0$.  Let $l(\lambda)$, the length of $\lambda$, be the number of nonzero $\lambda_{i}$ and let $|\lambda|$, the weight of $\lambda$, be the sum of the nonzero $\lambda_{i}$.  We will write $\lambda = \mu^{2}$ if there exists a partition $\mu$ such that $\lambda_{2i-1} = \lambda_{2i} = \mu_{i}$  (equivalently all parts of $\lambda$ occur with even multiplicity).  Analogously, we write $\lambda = 2\mu$ if there exists a partition $\mu$ such that $\lambda_{i} = 2\mu_{i}$ (equivalently each part of $\lambda$ is even).  Also let $m_{i}(\lambda)$ be the number of $\lambda_{j}$ equal to $i$ for each $i \geq 0$.

Recall the $t$-integer $[i] = [i]_{t} = (1-t^{i})/(1-t)$, as well as the $t$-factorial $[m]! = [m][m-1] \cdots [1]$, $[0]! = 1$.  Let 
\begin{align*}
\phi_{r}(t) &= (1-t)(1-t^{2}) \cdots (1-t^{r}),
\end{align*}
so that in particular $\phi_{r}(t)/(1-t)^{r} = [r]!$.  Then we define
\begin{align*}
v_{\lambda}(t) &= \prod_{i \geq 0} \prod_{j=1}^{m_{i}(\lambda)} \frac{1- t^{j}}{1-t}  
= \prod_{i \geq 0} \frac{\phi_{m_{i}(\lambda)}(t)}{(1-t)^{m_{i}(\lambda)}} = \prod_{i \geq 0} [m_{i}(\lambda)]!
\end{align*}
and
\begin{align*}
v_{\lambda+}(t) &= \prod_{i \geq 1} \prod_{j=1}^{m_{i}(\lambda)} \frac{1-t^{j}}{1-t} 
= \prod_{i \geq 1} \frac{\phi_{m_{i}(\lambda)}(t)}{(1-t)^{m_{i}(\lambda)}} = \prod_{i \geq 1} [m_{i}(\lambda)]!,
\end{align*}
so that the first takes into account the zero parts, while the second does not.  The Hall--Littlewood polynomial $P_{\lambda}(x_{1}, \dots, x_{n};t)$ indexed by $\lambda$ is defined to be
\begin{align*}
\frac{1}{v_{\lambda}(t)} \sum_{w \in S_{n}} w\Big( x^{\lambda} \prod_{1 \leq i<j \leq n} \frac{x_{i}-tx_{j}}{x_{i}-x_{j}} \Big),
\end{align*}
where we write $x^{\lambda}$ for $x_{1}^{\lambda_{1}} \cdots x_{n}^{\lambda_{n}}$ and $w$ acts on the subscripts of the $x_{i}$.  The normalization $1/v_\lambda(t)$ has the effect of making the 
coefficient of $x^\lambda$ equal to unity. (We will also write $P_{\lambda}^{(n)}(x;t)$ and use $P_{\lambda}(x^{(m)}, y^{(n)};t)$ to denote $P_{\lambda}(x_{1}, \dots, x_{m}, y_{1}, \dots, y_{n};t)$ in the final section).  We define the polynomials $\{R_{\lambda}^{(n)}(x;t)\}$ by $R_{\lambda}^{(n)}(x;t)=v_{\lambda}(t)P_{\lambda}^{(n)}(x;t)$.  For $w \in S_{n}$, we also define
\begin{equation} \label{Rpol}
R_{\lambda, w}^{(n)}(x;t) = w\Big( x^{\lambda} \prod_{1 \leq i<j \leq n} \frac{x_{i}-tx_{j}}{x_{i}-x_{j}} \Big),
\end{equation}
so that $R_{\lambda, w}^{(n)}(x;t)$ is the term of $R_{\lambda}^{(n)}(x;t)$ associated to the permutation $w$.  

There are two important degenerations of the Hall--Littlewood symmetric functions: at $t=0$, we recover the Schur functions $s_{\lambda}(x)$ and at $t=1$ the monomial symmetric functions $m_{\lambda}(x)$.  We remark that the Macdonald polynomials $P_{\lambda}(x;q,t)$ do not have poles at $q=0$, so there is no obstruction to specializing $q$ to zero; in fact we obtain the Hall--Littlewood polynomials (see \cite{Mac}, Ch. 6).  Similarly, when $q=t$ (or $q=0$ then $t=0$), $P_{\lambda}(x;q,t)$ reduces to $s_{\lambda}(x)$.  

Let 
\begin{align*}
b_{\lambda}(t) = \prod_{i \geq 1} \phi_{m_{i}(\lambda)}(t) 
= v_{\lambda+}(t) (1-t)^{l(\lambda)}.
\end{align*}
Then we let $Q_{\lambda}(x;t)$ be multiples of the $P_{\lambda}(x;t)$:
\begin{align*}
Q_{\lambda}(x;t) &= b_{\lambda}(t) P_{\lambda}(x;t);
\end{align*}
these form the adjoint basis with respect to the $t$-analog of the Hall inner product.  With this notation the Cauchy identity for Hall--Littlewood functions is
\begin{equation} \label{Cauchyid}
\sum_{\lambda} P_{\lambda}(x;t) Q_{\lambda}(x;t) = \prod_{i,j \geq 1} \frac{1-tx_{i}y_{j}}{1-x_{i}y_{j}}.  
\end{equation}

We recall the definition of Rogers--Szeg\H{o} polynomials, which appears in Sections 5--7.  Let $m$ be a nonnegative integer.  Then we let $H_{m}(z;t)$ denote the Rogers--Szeg\H{o} polynomial (see \cite{A}, Ch. 3, Examples 3--9)
\begin{equation}  \label{RSpol}
H_{m}(z;t) = \sum_{i=0}^{m} z^{i} {\qbinom{m}{i}}_{t},
\end{equation}
where 
\begin{equation*}
{\qbinom{m}{i}}_{t} = \begin{cases}
\frac{[m]!}{[m-i]![i]!}, & \text{ if } m \geq i \geq 0 \\
0, &\text{ otherwise}
\end{cases}
\end{equation*}
is the $t$-binomial coefficient.  It can be verified that the Rogers--Szeg\H{o} polynomials satisfy the following second-order recurrence:
\begin{align*}
H_{m}(z;t) &= (1+z)H_{m-1}(z;t) - (1-t^{m-1})zH_{m-2}(z;t).
\end{align*}

Also, we define the symmetric $q=0$ Selberg density \cite{RV}:
\begin{align*}
\tilde \Delta_{S}^{(n)}(x;t) &= \prod_{1 \leq i \neq j \leq n} \frac{1-x_{i}x_{j}^{-1}}{1-tx_{i}x_{j}^{-1}}
\end{align*}
and the symmetric Koornwinder density \cite{K}:
\begin{equation} \label{Kd}
\tilde \Delta_{K}^{(n)}(x;a,b,c,d;t) = \frac{1}{2^{n}n!} \prod_{1 \leq i \leq n} \frac{1-x_{i}^{\pm 2}}{(1-ax_{i}^{\pm 1})(1-bx_{i}^{\pm 1})(1-cx_{i}^{\pm 1})(1-dx_{i}^{\pm 1})} \prod_{1 \leq i<j \leq n} \frac{1-x_{i}^{\pm 1}x_{j}^{\pm 1}}{1-tx_{i}^{\pm 1}x_{j}^{\pm 1}},
\end{equation}
where we write $1-x_{i}^{\pm 2}$ for the product $(1-x_{i}^{2})(1-x_{i}^{-2})$ and $1-x_{i}^{\pm 1}x_{j}^{\pm 1}$ for $(1-x_{i}x_{j})(1-x_{i}^{-1}x_{j}^{-1})(1-x_{i}^{-1}x_{j})(1-x_{i}x_{j}^{-1})$ etc.  For convenience, we will write $\tilde \Delta_{S}^{(n)}$ and $\tilde \Delta_{K}^{(n)}(a,b,c,d)$ with the assumption that these densities are in $x_{1}, \dots, x_{n}$ with parameter $t$ when it is clear.    
We recall some notation for hypergeometric series from \cite{RV} and \cite{R}.  We define the $q$-symbol
\begin{align*}
(a;q) &= \prod_{k \geq 0} (1-aq^{k})
\end{align*}
and $(a_{1}, a_{2}, \dots, a_{l};q) = (a_{1};q)(a_{2};q) \cdots (a_{l};q)$.  Also, let 
\begin{equation*}
(a;q)_{n} = \prod_{j=0}^{n-1} (1-aq^{j})
\end{equation*}
for $n>0$ and $(a;q)_{0} = 1$.  We also define the $C$-symbols, which appear in the identities of \cite{RV}.  Let
\begin{align*}
C^{0}_{\mu}(x;q,t) &= \prod_{1 \leq i \leq l(\mu)} \frac{(t^{1-i}x;q)}{(q^{\mu_{i}}t^{1-i}x;q)} \\
C^{-}_{\mu}(x;q,t) &= \prod_{1 \leq i \leq l(\mu)} \frac{(x;q)}{(q^{\mu_{i}}t^{l(\mu)-i}x;q)}  \prod_{1 \leq i<j \leq l(\mu)} \frac{(q^{\mu_{i}-\mu_{j}}t^{j-i}x;q)}{(q^{\mu_{i}-\mu_{j}}t^{j-i-1}x;q)} \\
C^{+}_{\mu}(x;q,t) &= \prod_{1 \leq i \leq l(\mu)} \frac{(q^{\mu_{i}}t^{2-l(\mu)-i}x;q)}{(q^{2\mu_{i}}t^{2-2i}x;q)} \prod_{1 \leq i<j \leq l(\mu)} \frac{(q^{\mu_{i}+\mu_{j}}t^{3-j-i}x;q)}{(q^{\mu_{i}+\mu_{j}}t^{2-j-i}x;q)}.
\end{align*}
We note that $C_{\mu}^{0}(x;q,t)$ is the $q,t$-shifted factorial.  As before, we extend this by $C^{0, \pm}_{\mu}(a_{1}, a_{2}, \dots, a_{l};q) = C^{0, \pm}_{\mu}(a_{1};q) \cdots C^{0, \pm}_{\mu}(a_{l};q)$.  

We note that for $q=0$ we have
\begin{align*}
C^{0}_{\mu}(x;0,t) &= \prod_{1 \leq i \leq l(\mu)} (1-t^{1-i}x) \\
C^{-}_{\mu}(t;0,t) &= (1-t)^{l(\mu)}v_{\mu+}(t)\\
C^{+}_{\mu}(x;0,t) &=  1.\\
\end{align*}
Finally, we explain some notation involving permutations.  Let $w \in S_{n}$ act on the variables $z_{1}, \dots, z_{n}$ by 
\begin{align*}
w(z_{1} \cdots z_{n}) = z_{w(1)} \cdots z_{w(n)}
\end{align*}
as in the definition of Hall--Littlewood polynomials above.  We view the permutation $w$ as this string of variables.  For example the condition ``$z_{i}$ is in the $k$-th position of $w$" means that $w(k) = i$.  Also we write 
\begin{align*}
``z_{i} \prec_{w} z_{j}"
\end{align*}
if $i = w(i')$ and $j= w(j')$ for some $i' < j'$, i.e., $z_{i}$ appears to the left of $z_{j}$ in the permutation representation $z_{w(1)} \cdots z_{w(n)}$.  For $w \in S_{2n}$, we use $w(x_{1}^{\pm 1}, \dots, x_{n}^{\pm 1})$ to represent $z_{w(1)} \cdots z_{w(2n)}$, with $z_{i} = x_{i}$ for $1 \leq i \leq n$ and $z_{j} = x_{j-n}^{-1}$ for $n+1 \leq j \leq 2n$.

\section{Hall--Littlewood Orthogonality}
It is a well known result that Hall--Littlewood polynomials are orthogonal with respect to the density $\tilde \Delta_{S}$.  We prove this result using our method below, to illustrate the technique in a simple case.  

\begin{theorem}\label{orthog}
We have the following orthogonality relation for Hall--Littlewood polynomials:
\begin{equation*}
\int_{T} P_{\lambda}(x_{1}, \dots, x_{n};t) P_{\mu}(x_{1}^{-1}, \dots, x_{n}^{-1};t) \tilde \Delta_{S}^{(n)}(x;t) dT = \delta_{\lambda \mu} \frac{n!}{v_{\mu}(t)} 
\end{equation*}
\end{theorem}

\begin{proof}
Note first that by the definition of Hall--Littlewood polynomials, the LHS is a sum of $(n!)^{2}$ integrals in bijection with $S_{n} \times S_{n}$.  Now, since the integral is invariant under inverting all variables, we may restrict to the case where $\lambda \geq \mu$ in the reverse lexicographic ordering (we assume this throughout).  We will show that each of these terms vanish unless $\lambda = \mu$, and this argument will allow us to compute the normalization in the case $\lambda = \mu$.  By symmetry and (\ref{Rpol}), we have
\begin{equation*}
\int_{T} P_{\lambda}^{(n)}(x;t) P_{\mu}^{(n)}(x^{-1};t) \tilde \Delta_{S}^{(n)} dT = \frac{n!}{v_{\lambda}(t)v_{\mu}(t)} \sum_{\rho \in S_{n}} \int_{T} R_{\lambda, \text{id}}^{(n)}(x;t) R_{\mu, \rho}^{(n)}(x^{-1};t) \tilde \Delta_{S}^{(n)} dT.
\end{equation*}

\begin{claim}
We have the term-evaluation
\begin{equation*}
 \int_{T} R_{\lambda, \text{id}}^{(n)}(x;t) R_{\mu, \rho}^{(n)}(x^{-1};t) \tilde \Delta_{S}^{(n)} dT = 
t^{i(\rho)} 
\end{equation*}
if $x_{1}^{\lambda_{1}} \cdots x_{n}^{\lambda_{n}}x_{\rho(1)}^{-\mu_{1}} \cdots x_{\rho(n)}^{-\mu_{n}} =1$, and is otherwise equal to 0.   Here $i(\rho)$ is the number of inversions of $\rho$ with respect to the permutation $x_{1}^{-1} \cdots x_{n}^{-1}$.
\end{claim}
\noindent Note that $i(\rho)$ is the Coxeter length and recall the distribution of this statistic: $\sum_{\rho} t^{i(\rho)} = [n]!$.

To prove the claim, we use induction on $n$.  Note first that for $n=1$, the only term is $\int x_{1}^{\lambda_{1}} x_{1}^{-\mu_{1}} dT$, which vanishes unless $\lambda_{1} = \mu_{1}$.  Now suppose the result is true for $n-1$.  With this assumption we want to show that it holds true for $n$ variables. One can compute, by integrating with respect to $x_{1}$ in the iterated integral, that the LHS above is equal to
\begin{equation*}
\int_{T_{n-1}} \Big( \int_{T_{1}} x_{1}^{\lambda_{1}-\mu_{\rho^{-1}(1)}} \prod_{x_{j}^{-1} \prec_{\rho} x_{1}^{-1}} \frac{tx_{j} - x_{1}}{x_{j} - tx_{1}}\frac{dx_{1}}{2\pi \sqrt{-1}x_{1}}\Big) R_{\widehat{\lambda}, \widehat{\text{id}}}^{(n-1)}(x;t) R_{\widehat{\mu}, \widehat{\rho}}^{(n-1)}(x^{-1};t) \tilde \Delta_{S}^{(n-1)}(x;t) dT, 
\end{equation*}
where
\begin{align*}
\widehat{\text{id}} &= \text{id} \text{ with $x_{1}$ deleted} \\
\widehat{\rho} &= \rho \text{ with $x_{1}^{-1}$ deleted} \\
\widehat{\lambda} &= \lambda \text{ with $\lambda_{1}$ deleted} \\
\widehat{\mu} &= \mu \text{ with $\mu_{\rho^{-1}(1)}$ deleted}.
\end{align*}
Recall that $\lambda_{1} \geq \mu_{1} \geq \mu_{i}$ for all $1 \leq i \leq n$.  Thus, the inner integral in $x_{1}$ is zero if $\lambda_{1} > \mu_{\rho^{-1}(1)}$ and is $t^{|\{ j: x_{j}^{-1} \prec_{\rho} x_{1}^{-1} \}|}$ if $\lambda_{1} = \mu_{\rho^{-1}(1)}$.  In the latter case, note that $\widehat{\lambda} \geq \widehat{\mu}$, so we may use the induction hypothesis on the resulting $(n-1)$-dimensional integral, and combining this with the contribution from $x_{1}$ gives the result of the claim.

Note that the claim implies each term is zero if $\lambda \neq \mu$, so consequently the entire integral in zero.  Finally, we use the claim to compute the normalization value in the case $\lambda = \mu$.  By the above remarks, we have
\begin{equation*}
\int_{T} P_{\lambda}^{(n)}(x;t) P_{\mu}^{(n)}(x^{-1};t) \tilde \Delta_{S}^{(n)} dT = \frac{n!}{v_{\mu}(t)^{2}}\sum_{\substack{\rho \in S_{n} : \\ x_{1}^{\lambda_{1}} \cdots x_{n}^{\lambda_{n}}x_{\rho(1)}^{-\mu_{1}} \cdots x_{\rho(n)}^{-\mu_{n}}=1}} t^{i(\rho)}
\end{equation*}
Note that the permutations in the index of the sum are in statistic-preserving bijection with $S_{m_{0}(\mu)} \times S_{m_{1}(\mu)} \times \cdots$ so, using the comment immediately following the Claim, the above expression is equal to
\begin{equation*}
\frac{n!}{v_{\mu}(t)^{2}} \sum_{\rho \in S_{m_{0}(\mu)} \times S_{m_{1}(\mu)} \times \cdots} t^{i(\rho)} = \frac{n!}{v_{\mu}(t)^{2}} \prod_{i \geq 0} [m_{i}(\mu)]! = \frac{n!}{v_{\mu}(t)},
\end{equation*}
as desired.
\end{proof}

\section{$\alpha$ version}
In this section, we prove the orthogonal group integrals with an extra parameter $\alpha$.  This gives four identities - one for each component of $O(l)$, depending on the parity of $l$.  First, we use a result of Gustafson \cite{G} to compute some normalizations that will be used throughout the paper.  
\begin{proposition} \label{normalizations}
We have the following normalizations:
\begin{enumerate}
\item (\text{symplectic})
\begin{align*}
\int_{T} \tilde \Delta_{K}^{(n)}(x;\pm \sqrt{t},0,0;t) dT &= \frac{(1-t)^{n}}{(t^{2};t^{2})_{n}}
\end{align*}
\item (\text{Kawanaka})
\begin{align*}
 \int_{T} \tilde \Delta_{K}^{(n)}(x;1, \sqrt{t}, 0,0;t) dT &=  \frac{(1-t)^{n}}{(\sqrt{t};\sqrt{t})_{2n}}
 \end{align*}
 \item ($O^{+}(2n)$)
\begin{align*}
 \int_{T} \tilde \Delta_{K}^{(n)}(x;\pm 1, \pm \sqrt{t};t) dT &= \frac{(1-t)^{n}}{2(t;t)_{2n}}
 \end{align*}
 \item ($O^{-}(2n)$)
\begin{align*}
 \int_{T} \tilde \Delta_{K}^{(n-1)}(x;\pm t, \pm \sqrt{t};t) dT &= \frac{(1-t)^{n-1}}{(t^{3};t)_{2n-2}}
 \end{align*}
 \item ($O^{+}(2n+1)$)
\begin{align*}
\int_{T} \tilde \Delta_{K}^{(n)}(x;t,-1, \pm \sqrt{t};t) dT &= \frac{(1-t)^{n+1}}{(t;t)_{2n+1}}
\end{align*}
\item ($O^{-}(2n+1)$)
\begin{align*}
\int_{T} \tilde \Delta_{K}^{(n)}(x;1,-t, \pm \sqrt{t};t) dT &= \frac{(1-t)^{n+1}}{(t;t)_{2n+1}}.
\end{align*}
\end{enumerate}
\end{proposition}
We omit the proof, but in all cases it follows from setting $q=0$ and the appropriate values of $(a,b,c,d)$ in the integral evaluation:
\begin{align*}
\int_{T} \tilde \Delta_{K}^{(n)}(x;a,b,c,d;q,t) dT &= \prod_{0 \leq j<n} \frac{(t,t^{2n-2-j}abcd;q)}{(t^{j+1},t^{j}ab, t^{j}ac, t^{j}ad, t^{j}bc, t^{j}bd, t^{j}cd;q)},
\end{align*}
which may be found in \cite{G}.  

We remark that at $t=0$ the above densities have special significance.  In particular, (i) is the eigenvalue density of the symplectic group and (iii) - (vi) are the eigenvalue densities of $O^{+}(2n), O^{-}(2n), O^{+}(2n+1)$ and $O^{-}(2n+1)$ (in the orthogonal group case, the density depends on the component of the orthogonal group as well as whether the dimension is odd or even).  The density in (ii) appears in Corollary \ref{Kawanaka}, and that result corresponds to a summation identity of Kawanaka \cite{Ka1}.

In this section, we want to use a technique similar to the one used to prove Hall--Littlewood orthogonality.  Namely, we want to break up the integral into a sum of terms, one for each permutation, and study the resulting term integral.  The obstruction to this approach is that in many cases the poles lie on the contour, i.e., occur at $\pm 1$, so the pieces of the integral are not well-defined.  However, since the overall integral does not have singularities, we may use the principal value integral which we denote by P.V. (see \cite{Kan}, Section 8.3).  We first prove some results involving the principal value integrals.

\begin{lemma} \label{pvl1}
 Let $f(z)$ be a function in $z$ such that $zf(z)$ is holomorphic in a neighborhood of the unit disk.  Then
$$
\PV \int_{T} f(z) \frac{1}{1-z^{-2}} dT = \frac{f(1) + f(-1)}{4}.
$$
\end{lemma}

\begin{proof}
We have
\begin{multline*}
\PV \frac{1}{2\pi \sqrt{-1}} \int_{|z| =1} f(z) \frac{1}{1-z^{-2}} \frac{1}{z} dz
 = \lim_{\epsilon \rightarrow 0^{+}} \frac{1}{2} \Biggl[ \frac{1}{2\pi \sqrt{-1}} \int_{|z| = 1-\epsilon} zf(z) \frac{1}{z^{2} - 1}dz \\
 + \frac{1}{2\pi \sqrt{-1}} \int_{|z| = 1 + \epsilon} zf(z) \frac{1}{z^{2} - 1} dz \Biggr]
\end{multline*}
But now as $zf(z)$ is holomorphic in a neighborhood of the disk, and the singularities of $1/(z^{2}-1)$ lie outside of the disk, the first integral is zero by Cauchy's theorem.  Using the residue theorem for the second integral (it has simple poles at $\pm 1$) gives
\begin{align*}
\lim_{\epsilon \rightarrow 0} \frac{1}{2} \Biggl[ \text{Res}_{z=1} \frac{zf(z)}{(z-1)(z+1)} + \text{Res}_{z=-1} \frac{zf(z)}{(z-1)(z+1)} \Biggr] 
= \frac{1}{2}\Biggl[\frac{f(1)}{2} + \frac{f(-1)}{2}\Biggr] 
= \frac{1}{4}\Big[f(1) + f(-1)\Big].
\end{align*}
\end{proof}

\begin{lemma} \label{pvl2}
Let $p$ be a function in $x_{1}, \dots, x_{n}$ such that $x_{i}p$ is holomorphic in $x_{i}$ in a neighborhood of the unit disk for all $1 \leq i \leq n$ and $p( \pm 1, \dots, \pm 1) = 0$ for all $2^{n}$ combinations.  Let $\Delta$ be a function in $x_{1}, \dots, x_{n}$ such that $\Delta( \pm 1, \dots, \pm 1, x_{i+1}, \dots, x_{n})$ is holomorphic in $x_{i+1}$ in a neighborhood of the unit disk for all $0 \leq i \leq n-1$ (again for all $2^{i}$ combinations).  Then
$$
\PV \int_{T} p \cdot \Delta \cdot \prod_{1 \leq i \leq n} \frac{1}{1-x_{i}^{-2}} dT = 0.
$$  
\end{lemma}

\begin{proof}
We give a proof by induction on $n$.
For $n=1$, since $x_{1} \cdot p \cdot \Delta$ is holomorphic in $x_{1}$ we may use Lemma \ref{pvl1}:
$$
\PV \int_{T} p \cdot \Delta \cdot \frac{1}{1-x_{1}^{-2}} dT = \frac{1}{4}[ p(1)\Delta(1) + p(-1)\Delta(-1)].
$$
But then $p(1) = p(-1) = 0$ by assumption, so the integral is zero as desired.  

Now suppose the result holds in the case of $n-1$ variables.  Consider the $n$ variable case, and let $p, \Delta$ in $x_{1}, \dots, x_{n}$ satisfy the above conditions.  Integrate first with respect to $x_{1}$ and note that $x_{1} \cdot p \cdot \Delta$ is holomorphic in $x_{1}$ so we can apply Lemma \ref{pvl1}:
\begin{multline*}
\PV \int_{T} p \cdot \Delta \cdot \prod_{1 \leq i \leq n} \frac{1}{1-x_{i}^{-2}} dT = \frac{1}{4} \PV \int_{T_{n-1}} p(1,x_{2}, \dots, x_{n})\Delta(1,x_{2}, \dots,x_{n}) \prod_{2 \leq i \leq n} \frac{1}{1-x_{i}^{-2}} dT \\
+ \frac{1}{4} \PV \int_{T_{n-1}} p(-1, x_{2}, \dots, x_{n})\Delta(-1, x_{2}, \dots, x_{n}) \prod_{2 \leq i \leq n} \frac{1}{1-x_{i}^{-2}} dT.
\end{multline*}
But now the pairs $p(1,x_{2}, \dots, x_{n}), \Delta(1,x_{2},\dots,x_{n})$ and $p(-1,x_{2},\dots,x_{n}),$
$\Delta(-1,x_{2},\dots,x_{n})$ satisfy the conditions of the theorem for $n-1$ variables $x_{2}, \dots, x_{n}$, so by the induction hypothesis each of the two integrals is zero, so the total integral is zero. 
\end{proof}

For this section, we let $\rho_{2n} = (1, 2, \dots, 2n)$.  We also let $1^{k} = (1,1, \dots, 1)$ with exactly $k$ ones.  As above we will work with principal value integrals, as necessary.  For simplicity, we will suppress the notation P.V.

\begin{theorem}\label{ape}
Let $l(\lambda) \leq 2n$.  We have the following integral identity for $O^{+}(2n)$:
\begin{multline*}
\frac{1}{\int \tilde \Delta_{K}^{(n)}(\pm 1, \pm \sqrt{t}) dT} \int P_{\lambda}(x_{1}^{\pm 1}, \dots, x_{n}^{\pm 1};t) \tilde \Delta_{K}^{(n)}( \pm 1, \pm \sqrt{t}) \prod_{i=1}^{n} (1-\alpha x_{i}^{\pm 1})dT\\
= \frac{\phi_{2n}(t)}{v_{\lambda}(t) (1-t)^{2n}} \Big[ (-\alpha)^{\text{\# of odd parts of $\lambda$}} + (-\alpha)^{\# \text{ of even parts of $\lambda$}}\Big] \\
= \frac{[2n]!}{v_{\lambda}(t)} \Big[ (-\alpha)^{\text{\# of odd parts of $\lambda$}} + (-\alpha)^{\# \text{ of even parts of $\lambda$}}\Big].
\end{multline*}
\end{theorem}

\begin{proof}
We will first show the following:
\begin{align*}
\int R_{\lambda}(x_{1}^{\pm 1}, \dots, x_{n}^{\pm 1};t) \tilde \Delta_{K}^{(n)}( \pm 1, \pm \sqrt{t}) \prod_{i=1}^{n} (1-\alpha x_{i}^{\pm 1})dT
= \frac{1}{2^{n}(1-t)^{n}}\mathrm{Pf}[a_{j,k}]^{\lambda},
\end{align*}
where $\mathrm{Pf}$ denotes the Pfaffian and the $2n \times 2n$ antisymmetric matrix $[a_{j,k}]^{\lambda}$ is defined by
\begin{align*}
a_{j,k}^{\lambda} &= (1+\alpha^{2})\chi_{(\lambda_{j}-j)-(\lambda_{k}-k) \text{ odd}} + 2(-\alpha)\chi_{(\lambda_{j}-j)-(\lambda_{k}-k) \text{ even}},
\end{align*}
for $1 \leq j<k \leq 2n$.

First, note that by symmetry we can rewrite the above integral as $2^{n} n!$ times the sum over all matchings $w$ of $x_{1}^{\pm 1}, \dots, x_{n}^{\pm 1}$, where a matching is a permutation in $S_{2n}$ such that $x_{i}$ occurs to the left of $x_{i}^{-1}$ and $x_{i}$ occurs to the left of $x_{j}$ for $1 \leq i <j \leq n$.  In particular, $x_{1}$ occurs first.  Thus, we have 
\begin{multline*}
 \int R_{\lambda}^{(2n)}(x^{\pm 1};t) \tilde \Delta_{K}^{(n)}(\pm 1; \pm \sqrt{t}) \prod_{i=1}^{n}(1-\alpha x_{i}^{\pm 1})dT\\
= 2^{n}n! \sum_{w} \int R_{\lambda, w}^{(2n)}(x^{\pm 1};t) \tilde \Delta_{K}^{(n)}(\pm 1; \pm \sqrt{t}) \prod_{i=1}^{n}(1-\alpha x_{i}^{\pm 1})dT,
\end{multline*}
where the sum is over matchings $w$ in $S_{2n}$.  

We introduce some notation for a matching $w \in S_{2n}$.  We write $w = \{ (i_{1}, i_{1}'), \dots, (i_{n}, i_{n}') \}$ to indicate that $x_{k}$ occurs in position $i_{k}$ and $x_{k}^{-1}$ occurs in position $i_{k}'$ for all $1 \leq k \leq n$.  Clearly we have $i_{k} < i_{k}'$ for all $k$ and $i_{j} < i_{k}$ for all $j<k$.  

\begin{claim}\label{apet}
Let $\lambda = (\lambda_{1}, \dots, \lambda_{2n})$ with $\lambda_{1} \geq \lambda_{2} \geq \cdots \geq \lambda_{2n} \in \mathbb{Z}$.  Then we have the following term-evaluation:
\begin{align*}
2^{n}n! \PV \int_{T} R_{\lambda, w}(x_{1}^{\pm 1}, \dots, x_{n}^{\pm 1};t) \tilde \Delta_{K}^{(n)}(\pm 1, \pm \sqrt{t}) \prod_{i=1}^{n} (1-\alpha x_{i}^{\pm 1}) dT &= \frac{\epsilon(w)}{2^{n}(1-t)^{n}} \prod_{1 \leq k \leq n} a_{i_{k},i_{k}'}^{\lambda},
\end{align*}
where 
$\epsilon(w)$ is the sign of $w$ and $a_{i_{k},i_{k}'}^{\lambda}$ is the $(i_{k}, i_{k}')$ entry of the matrix $[a_{j,k}]^{\lambda}$.  In particular, the term integral only depends on the parity of the parts $\lambda_{1}, \dots, \lambda_{2n}$.  
\end{claim}

Let $\mu$ be such that $\lambda = \mu + \rho_{2n}$.  We give a proof by induction on $n$, the number of variables.  For $n=1$, there is only one matching---in particular, $x_{1}^{-1}$ must occur in position $2$. The (principal value) integral is
\begin{multline*}
\int_{T} x_{1}^{\lambda_{1} - \lambda_{2}} \frac{(1-tx_{1}^{-2})}{(1-x_{1}^{-2})} \frac{(1-\alpha x_{1})(1-\alpha x_{1}^{-1})}{(1-tx_{1}^{2})(1-tx_{1}^{-2})} dT
= \int_{T} x_{1}^{\lambda_{1} - \lambda_{2}} \frac{(1-\alpha x_{1})(1-\alpha x_{1}^{-1})}{(1-x_{1}^{-2})(1-tx_{1}^{2})} dT \\
= \int_{T} x_{1}^{\lambda_{1}-\lambda_{2}} \frac{(1 + \alpha^{2}) - \alpha(x_{1} + x_{1}^{-1})}{(1-tx_{1}^{2})(1-x_{1}^{-2})} dT
\end{multline*}
and $\lambda_{1}-\lambda_{2} \geq 0$.  Note that the conditions for Lemma \ref{pvl1} are satisfied.  Applying that result gives that the value of the integral is $2(-\alpha)/2(1-t)$ if $\lambda_{1}-\lambda_{2}$ is odd, and $(1+\alpha^{2})/2(1-t)$ if $\lambda_{1}-\lambda_{2}$ is even, which agrees with the above claim.  

Now suppose the result is true for up to $n-1$ variables and consider the $n$ variable case.  Note first that $i_{1} = 1$.  One can compute, by combining terms involving $x_{1}$ in the iterated integral, that
\begin{multline*}
2^{n}n! \int R_{\lambda, w}^{(2n)}(x^{\pm 1};t) \tilde \Delta_{K}^{(n)}(\pm 1, \pm \sqrt{t}) \prod_{i=1}^{n} (1-\alpha x_{i}^{\pm 1}) dT \\
=\int_{T_{n-1}} \Big( \int_{T_{1}} x_{1}^{\lambda_{1}-\lambda_{i_{1}'}}  \frac{(1-\alpha x_{1})(1-\alpha x_{1}^{-1})}{(1-tx_{1}^{2})(1-x_{1}^{-2})} \prod_{\substack{x_{j} : \\ x_{1} \prec_{w} x_{j} \prec_{w} x_{1}^{-1} \prec_{w} x_{j}^{-1} }} \frac{(t-x_{1}x_{j})}{(1-tx_{1}x_{j})}\\
 \prod_{\substack{x_{j} : \\ x_{1} \prec_{w} x_{j} \prec_{w} x_{j}^{-1} \prec_{w} x_{1}^{-1}}}\frac{(t-x_{1}x_{j}^{-1})(t-x_{1}x_{j})}{(1-tx_{1}x_{j}^{-1})(1-tx_{1}x_{j})} dT \Big) F_{\widehat{\lambda}, \tilde w} dT,
\end{multline*}
where 
\begin{equation*}
F_{\widehat \lambda, \tilde w} = 2^{n-1}(n-1)! R_{\widehat \lambda, \tilde w}(x_{2}^{\pm 1}, \dots, x_{n}^{\pm 1};t) \tilde \Delta_{K}^{(n-1)}(\pm 1, \pm \sqrt{t}) \prod_{i=2}^{n} (1-\alpha x_{i}^{\pm 1})
\end{equation*}
and  $\widehat \lambda$ is $\lambda$ with parts $\lambda_{1}, \lambda_{i_{1}'}$ deleted; $\tilde w$ is $w$ with $x_{1}, x_{1}^{-1}$ deleted.  

In particular, the conditions for Lemma \ref{pvl1} are satisfied for the inner integral in $x_{1}$.  Note that the  terms 
$$
\frac{(t-x_{1}x_{i})}{(1-tx_{1}x_{i})}\frac{(t-x_{1}x_{i}^{-1})}{(1-tx_{1}x_{i}^{-1})}
$$
give $1$ when evaluated at $x_{1} = \pm 1$, so the above integral evaluates to
\begin{multline*}
\frac{1}{4(1-t)} \int_{T_{n-1}} \Bigg[ F_{\widehat \lambda, \tilde w} \cdot (1 + \alpha^{2} - 2\alpha) \Big(\prod_{\substack{x_{j} : \\ x_{1} \prec_{w} x_{j} \prec_{w} x_{1}^{-1} \prec_{w} x_{j}^{-1} }} \frac{t-x_{j}}{1-tx_{j}} \Big)\\
 + F_{\widehat \lambda, \tilde w} \cdot(1 + \alpha^{2} + 2\alpha)(-1)^{\lambda_{1}-\lambda_{i_{1}'}} \Big(\prod_{\substack{x_{j}: \\ x_{1} \prec_{w} x_{j} \prec_{w} x_{1}^{-1} \prec_{w} x_{j}^{-1}}}\frac{t+x_{j}}{1+tx_{j}} \Big) \Bigg] dT.
\end{multline*}
But now since $(t-x_{i})/(1-tx_{i})$ and $(t+x_{i})/(1+tx_{i})$ are power series in $x_{i}$, we may apply the inductive hypothesis to each part of the new integral: we reduce exponents on $x_{i}$ modulo $2$.  We get
\begin{multline*}
\frac{1}{4(1-t)} \int_{T_{n-1}} \Bigg[ F_{\widehat \lambda, \tilde w} \cdot (1 + \alpha^{2} - 2\alpha) \Big(\prod_{\substack{x_{j} : \\ x_{1} \prec_{w} x_{j} \prec_{w} x_{1}^{-1} \prec_{w} x_{j}^{-1} }} (-x_{j})\Big) \\
 + F_{\widehat \lambda, \tilde w} \cdot (1 + \alpha^{2} + 2\alpha) (-1)^{\lambda_{1}-\lambda_{i_{1}'}} \Big(\prod_{\substack{x_{j} : \\ x_{1} \prec_{w} x_{j} \prec_{w} x_{1}^{-1} \prec_{w} x_{j}^{-1} }} x_{j}\Big)\Bigg] dT.
\end{multline*}
But now note that 
\begin{equation*}
\prod_{\substack{x_{j} : \\ x_{1} \prec_{w} x_{j} \prec_{w} x_{1}^{-1} \prec_{w} x_{j}^{-1} }} (-1) = \prod_{\substack{x_{j} : \\ x_{1} \prec_{w} x_{j} \prec_{w} x_{1}^{-1} \prec_{w} x_{j}^{-1} }} (-1) \prod_{\substack{x_{j} : \\ x_{1} \prec_{w} x_{j} \prec_{w} x_{j}^{-1} \prec_{w} x_{1}^{-1} }} (-1)^{2} = (-1)^{i_{1}'-2},
\end{equation*}
since $i_{1}'-2$ is the number of variables between $x_{1}$ and $x_{1}^{-1}$ in the matching $w$.  We can compute
\begin{multline*}
(1 + \alpha^{2} - 2\alpha)(-1)^{i_{1}'-2}  + (1 + \alpha^{2} + 2\alpha)(-1)^{\lambda_{1}-\lambda_{i_{1}'}} = (1 + \alpha^{2})[(-1)^{i_{1}'} + (-1)^{\lambda_{1}-\lambda_{i_{1}'}} ]  - 2\alpha[(-1)^{i_{1}'} + (-1)^{\lambda_{1} -\lambda_{i_{1}'} +1}] \\ 
= \begin{cases} 2(-1)^{i_{1}'}(1+\alpha^{2}) & \text{if $\lambda_{1}-\lambda_{i_{1}'}+i_{1}'-1$ is odd,}
\\
 -4(-1)^{i_{1}'}\alpha &\text{if $\lambda_{1}-\lambda_{i_{1}'} +i_{1}'-1$ is even.} \\
\end{cases}
\end{multline*}
Combining this with the factor $1/4(1-t)$ and noting that 
\begin{align*}
F_{\widehat \lambda, \tilde w} \cdot \Big(\prod_{\substack{x_{j} : \\ x_{1} \prec_{w} x_{j} \prec_{w} x_{1}^{-1} \prec_{w} x_{j}^{-1} }} x_{j}\Big) &= F_{\tilde \lambda, \tilde w},
\end{align*}
with
\begin{align*}
\tilde \lambda &= (\lambda_{2}+1, \dots, \lambda_{i_{1}'-1}+1, \lambda_{i_{1}'+1}, \dots, \lambda_{2n}),
\end{align*}
gives that
\begin{multline*}
2^{n}n! \int R_{\lambda, w}^{(2n)}(x^{\pm 1};t) \tilde \Delta_{K}^{(n)}(\pm 1, \pm \sqrt{t}) \prod_{i=1}^{n} (1-\alpha x_{i}^{\pm 1}) dT \\
=  \frac{2^{n-1}(n-1)!}{2(1-t)} a_{i_{1}, i_{1}'}^{\lambda} (-1)^{i_{1}'}  \int_{T} R_{\tilde \lambda, \tilde w}(x_{1}^{\pm 1}, \dots, x_{n-1}^{\pm 1};t) \tilde \Delta_{K}^{(n-1)}(\pm 1, \pm \sqrt{t}) \prod_{i=1}^{n-1} (1-\alpha x_{i}^{\pm 1}) dT.
\end{multline*}
Now set $\widehat{\mu} = (\mu_{2}, \dots, \mu_{i_{1}'-1}, \mu_{i_{1}'+1}, \dots, \mu_{2n})$, and note that $\tilde \lambda$ and $\widehat{\mu} + \rho_{2n-2}$ have equivalent parts modulo $2$.  Thus, using the induction hypothesis twice, the above is equal to
\begin{multline*}
 \frac{2^{n-1}(n-1)!}{2(1-t)} a_{i_{1}, i_{1}'}^{\lambda} (-1)^{i_{1}'}  \int_{T} R_{\widehat{\mu} + \rho_{2n-2}, \tilde w}(x_{1}^{\pm 1}, \dots, x_{n-1}^{\pm 1};t) \tilde \Delta_{K}^{(n-1)}(\pm 1, \pm \sqrt{t}) \prod_{i=1}^{n-1} (1-\alpha x_{i}^{\pm 1}) dT \\
 =  \frac{a_{i_{1}, i_{1}'}^{\lambda} (-1)^{i_{1}'}}{2(1-t)} \frac{\epsilon(\tilde w)}{2^{n-1}(1-t)^{n-1}} \prod_{2 \leq k \leq n} a_{i_{k},i_{k}'}^{\widehat{\mu} + \rho_{2n-2}} = \frac{\epsilon(w)}{2^{n}(1-t)^{n}} \prod_{1 \leq k \leq n} a_{i_{k},i_{k}'}^{\lambda}
 \end{multline*}
 as desired.  This proves the claim.
 
Note in particular this result implies that the integral of a matching $w$ is the term in $\frac{1}{2^{n}(1-t)^{n}} \mathrm{Pf}[a_{j,k}]^{\lambda}$ corresponding to $w$.  

Now using the claim, we have
\begin{multline*}
 \int_{T} R_{\lambda}(x_{1}^{\pm 1}, \dots, x_{n}^{\pm 1};t) \tilde \Delta_{K}^{(n)}(\pm 1, \pm \sqrt{t}) \prod_{i=1}^{n} (1- \alpha x_{i}^{\pm 1}) dT \\
= 2^{n} n! \sum_{\substack{w \text{ a matching }\\ \text{in }S_{2n}}} \PV \int_{T} R_{\lambda, w}(x_{1}^{\pm 1}, \dots, x_{n}^{\pm 1};t) \tilde \Delta_{K}^{(n)}(\pm 1, \pm \sqrt{t}) \prod_{i=1}^{n} (1-\alpha x_{i}^{\pm 1}) dT 
= \frac{1}{2^{n}(1-t)^{n}} \mathrm{Pf}[a_{j,k}]^{\lambda}
\end{multline*}
since the term integrals are in bijection with the terms of the Pfaffian.  

Now we use this to prove the theorem.  Using Proposition \ref{normalizations}(iii), we have
\begin{multline*}
\frac{1}{\int \tilde \Delta_{K}^{(n)}(\pm 1, \pm \sqrt{t}) dT} \int P_{\lambda}(x_{1}^{\pm 1}, \dots, x_{n}^{\pm 1};t) \tilde \Delta_{K}^{(n)}( \pm 1, \pm \sqrt{t}) \prod_{i=1}^{n} (1-\alpha x_{i}^{\pm 1})dT\\
= \frac{2(1-t)(1-t^{2}) \cdots (1-t^{2n})}{(1-t)^{n}} \frac{1}{v_{\lambda}(t)2^{n}(1-t)^{n}}\mathrm{Pf}[a_{j,k}]^{\lambda} 
= \frac{(1-t)(1-t^{2}) \cdots (1-t^{2n})}{(1-t)^{2n}} \frac{1}{v_{\lambda}(t)2^{n-1}}\mathrm{Pf}[a_{j,k}]^{\lambda}.
\end{multline*}
But now by \cite[5.17]{FR}
\begin{align*}
\mathrm{Pf}[a_{j,k}]^{\lambda} &= 2^{n-1}\Big[ (-\alpha)^{\sum_{j=1}^{2n} [\lambda_{j} \text{ mod}2]}+ (-\alpha)^{\sum_{j=1}^{2n} [(\lambda_{j}+1) \text{ mod}2]}\Big],
\end{align*}
which gives the result.
\end{proof}

\begin{theorem}\label{ame}
Let $l(\lambda) \leq 2n$.  We have the following integral identity for $O^{-}(2n)$:
\begin{multline*}
\frac{(1-\alpha^{2})}{\int \tilde \Delta_{K}^{(n-1)}(\pm t, \pm \sqrt{t}) dT} \int P_{\lambda}(x_{1}^{\pm 1}, \dots, x_{n-1}^{\pm 1},1,-1;t) \tilde \Delta_{K}^{(n-1)}(\pm t, \pm \sqrt{t}) \prod_{i=1}^{n-1} (1-\alpha x_{i}^{\pm 1}) dT\\
= \frac{\phi_{2n}(t)}{v_{\lambda}(t) (1-t)^{2n}} \Big[ (-\alpha)^{\# \text{ of odd parts of $\lambda$}} - (-\alpha)^{\# \text{ of even parts of $\lambda$}}\Big]. \\
\end{multline*}
\end{theorem}

\begin{proof}
We will first show the following:
\begin{align*}
 \int R_{\lambda}(x_{1}^{\pm 1}, \dots, x_{n-1}^{\pm 1},1,-1;t) \tilde \Delta_{K}^{(n-1)}(\pm t, \pm \sqrt{t}) \prod_{i=1}^{n-1} (1-\alpha x_{i}^{\pm 1}) dT 
= \frac{(1+t)}{2} \frac{1}{2^{n-1}(1-t)^{n-1}}\mathrm{Pf}[M]^{\lambda},
\end{align*}
where the $(2n+2) \times (2n+2)$ antisymmetric matrix $[M]^{\lambda}$ is defined by 
\begin{align*}
\begin{cases} M_{1,2}^{\lambda} =0 \\
M_{1,k}^{\lambda} = (-1)^{\lambda_{k-2}-(k-2)} &\text{if $k \geq 3$} \\
M_{2,k}^{\lambda} = 1 &\text{if $k \geq 3$} \\
M_{j,k}^{\lambda} = a_{j-2,k-2}^{\lambda} &\text{if $3 \leq j<k \leq 2n+2$} 
\end{cases}
\end{align*}
and the $2n \times 2n$ matrix $[a_{j,k}]^{\lambda}$ is as in Theorem \ref{ape}.  

Note first that the integral is a sum of $(2n)!$ terms, but by symmetry we may restrict to the ``pseudo-matchings"---those with $\pm 1$ anywhere, but $x_{i}$ to the left of $x_{i}^{-1}$ for $1 \leq i \leq n-1$ and $x_{i}$ to the left of $x_{j}$ for $1 \leq i<j \leq n-1$.  There are $(2n)!/2^{n-1}(n-1)!$ such pseudo-matchings, and each has $2^{n-1}(n-1)!$ permutations with identical integral.  
\begin{claim} \label{ame1}
Let $w$ be a fixed pseudo-matching with $(-1)$ in position $j$ and $(+1)$ in position $k$ (here $1 \leq j \neq k \leq 2n$).  Then we have the following:
\begin{multline*}
 2^{n-1}(n-1)! \PV \int R_{\lambda, w}(x_{1}^{\pm 1}, \dots, x_{n-1}^{\pm 1}, \pm 1;t) \tilde \Delta_{K}^{(n-1)}(\pm t, \pm \sqrt{t}) \prod_{i=1}^{n-1} (1-\alpha x_{i}^{\pm 1}) dT \\
= 2^{n-1}(n-1)!  (-1)^{\lambda_{j}+k-2 + \chi_{j>k}} \frac{(1+t)}{2} \PV \int R_{\tilde \lambda, \tilde w}^{(2(n-1))}(x^{\pm 1};t) \tilde \Delta_{K}^{(n-1)}(\pm 1, \pm \sqrt{t}) \prod_{i=1}^{n-1}(1- \alpha x_{i}^{\pm 1})dT,
\end{multline*}
where $\tilde w$ is $w$ with $\pm 1$ deleted (in particular, a matching in $S_{2n-2}$) and $\tilde \lambda$ is $\lambda$ with parts $\lambda_{j}, \lambda_{k}$ deleted and all parts between $\lambda_{j}$ and $\lambda_{k}$ increased by $1$, so that (in the case $j<k$, for example)
\begin{align*}
\tilde \lambda &= (\lambda_{1}, \dots, \lambda_{j-1}, \lambda_{j+1}+1, \dots, \lambda_{k-1}+1, \lambda_{k+1}, \dots, \lambda_{2n}).
\end{align*}
\end{claim}
We prove the claim.  First, using (\ref{Kd}), we have
\begin{multline*}
2^{n-1}(n-1)!\tilde \Delta_{K}^{(n-1)}(\pm t, \pm \sqrt{t}) \\
=  \prod_{1 \leq i \leq n-1} \frac{1-x_{i}^{\pm 2}}{(1+tx_{i}^{\pm 1})(1-tx_{i}^{\pm 1})(1+\sqrt{t}x_{i}^{\pm 1})(1-\sqrt{t}x_{i}^{\pm 1})} \prod_{1 \leq i<j \leq n-1} \frac{1-x_{i}^{\pm 1}x_{j}^{\pm 1}}{1-tx_{i}^{\pm 1}x_{j}^{\pm 1}}.
\end{multline*}
Define the set $X = \{ (x_{i}^{\pm 1}, x_{j}^{\pm 1}) : 1 \leq i \neq j \leq n-1 \}$, and let $u_{\lambda, w}^{(n-1)}(x;t)$ be defined by
\begin{equation*}
R_{\lambda, w}(x_{1}^{\pm 1}, \dots, x_{n-1}^{\pm 1}, \pm 1;t) = u_{\lambda, w}^{(n-1)}(x;t) \prod_{\substack{(z_{i}, z_{j}) \in X: \\ z_{i} \prec_{w} z_{j}}} \frac{z_{i} - tz_{j}}{z_{i} - z_{j}}.  
\end{equation*}
Also define $p_{1}$ and $\Delta_{1}$ by
\begin{equation*}
 u_{\lambda, w}^{(n-1)}(x;t)  \prod_{1 \leq i \leq n-1} \frac{1-x_{i}^{\pm 2}}{(1+tx_{i}^{\pm 1})(1-tx_{i}^{\pm 1})(1+\sqrt{t}x_{i}^{\pm 1})(1-\sqrt{t}x_{i}^{\pm 1})}  \prod_{i=1}^{n-1} (1-\alpha x_{i}^{\pm 1})= p_{1} \prod_{i=1}^{n-1} \frac{1}{1-x_{i}^{-2}}
\end{equation*}
and 
\begin{equation*}
\prod_{1 \leq i<j \leq n-1} \frac{1-x_{i}^{\pm 1}x_{j}^{\pm 1}}{1-tx_{i}^{\pm 1}x_{j}^{\pm 1}}\prod_{\substack{(z_{i}, z_{j}) \in X: \\ z_{i} \prec_{w} z_{j}}} \frac{z_{i} - tz_{j}}{z_{i} - z_{j}} = \Delta_{1}.
\end{equation*}
Note that 
\begin{equation*}
R_{\lambda, w}(x_{1}^{\pm 1}, \dots, x_{n-1}^{\pm 1}, \pm 1;t) \tilde \Delta_{K}^{(n-1)}(\pm t, \pm \sqrt{t}) \prod_{i=1}^{n-1} (1-\alpha x_{i}^{\pm 1}) = p_{1} \Delta_{1} \prod_{i=1}^{n-1} \frac{1}{1-x_{i}^{-2}} .
\end{equation*}
Define analogously $p_{2}$ and $\Delta_{2}$ using $R_{\tilde \lambda, \tilde w}(x_{1}^{\pm 1}, \dots, x_{n-1}^{\pm 1};t)$ and $\tilde \Delta_{K}^{(n-1)}(\pm 1, \pm \sqrt{t})$ instead of $R_{\lambda, w}^{(2n)}(x^{\pm 1}, \pm 1;t)$ and $\tilde \Delta_{K}^{(n-1)}( \pm t, \pm \sqrt{t})$.  

Then one can check $\Delta_{1} = \Delta_{2} =: \Delta$ and $\Delta(\pm 1, \dots, \pm 1, x_{i+1}, \dots, x_{n-1})$ is holomorphic in $x_{i+1}$ for all $0 \leq i \leq n-2$ and all $2^{i}$ combinations.  Also, the function $p =p_{1} -  (-1)^{\lambda_{j}+k-2}\frac{(1+t)}{2}p_{2}$ (resp. $p=p_{1} - (-1)^{\lambda_{j}+k-1} \frac{(1+t)}{2}p_{2}$) satisfies the conditions of Lemma \ref{pvl2} if $j<k$ (resp. $j>k$).  So using that result, we have 
\begin{align*}
 \int p_{1} \cdot \Delta \cdot \prod_{1 \leq i \leq n-1} \frac{1}{1-x_{i}^{-2}} dT = (-1)^{\lambda_{j}+k-2} \frac{(1+t)}{2}\int p_{2} \cdot \Delta \cdot \prod_{1 \leq i \leq n-1} \frac{1}{1-x_{i}^{-2}} dT
\end{align*}
if $j<k$ and 
\begin{align*}
\int p_{1} \cdot \Delta \cdot \prod_{1 \leq i \leq n-1} \frac{1}{1-x_{i}^{-2}} dT = (-1)^{\lambda_{j}+k-1} \frac{(1+t)}{2} \int p_{2} \cdot \Delta \cdot \prod_{1 \leq i \leq n-1} \frac{1}{1-x_{i}^{-2}} dT
\end{align*}
if $j>k$.   Thus, in the case $j<k$ we obtain
\begin{multline*}
 \int R_{\lambda, w}(x_{1}^{\pm 1}, \dots, x_{n-1}^{\pm 1}, \pm 1;t) \tilde \Delta_{K}^{(n-1)}(\pm t, \pm \sqrt{t}) \prod_{i=1}^{n-1} (1-\alpha x_{i}^{\pm 1}) dT \\
= (-1)^{\lambda_{j}+k-2} \frac{(1+t)}{2}  \int R_{\tilde \lambda, \tilde w}(x_{1}^{\pm 1}, \dots, x_{n-1}^{\pm 1};t) \tilde \Delta_{K}^{(n-1)}(\pm 1, \pm \sqrt{t}) \prod_{i=1}^{n-1}(1- \alpha x_{i}^{\pm 1}) dT,
\end{multline*}
and analogously for the case $j>k$, which proves the claim.

As in Theorem \ref{ape}, we introduce notation for pseudo-matchings.  We write $\{(j,k), (i_{1}, i_{1}'), \dots, (i_{n-1}, i_{n-1}') \}$ for the pseudo-matching with $-1$ in position $j$, $1$ in position $k$ and $x_{k}$ in position $i_{k}$, $x_{k}^{-1}$ in position $i_{k}'$ for all $1 \leq k \leq n-1$.  Note that we have $i_{k} < i_{k}'$ and $i_{l} < i_{k}$ for $l<k$.  We may extend this to a matching in $S_{2(n+1)}$ by $\{(1,j+2), (2,k+2), (i_{1} + 2, i_{1}'+2), \dots, (i_{n-1}+2, i_{n-1}'+2) \} = \{ (j_{1}=1 , j_{1}' = j+2), (j_{2} =2, j_{2}' = k+2), \dots, (j_{n+1},j_{n+1}') \}$, with $i_{k}+2 = j_{k+2}$ and $i_{k}'+2 = j_{k+2}'$ for all $1 \leq k \leq n-1$.  

\begin{claim}\label{amet}
Let $w = \{ (j,k), (i_{1}, i_{1}'), \dots, (i_{n-1},i_{n-1}') \}$ be a pseudo-matching in $S_{2n}$, and extend it to a matching  $\{ (j_{1}=1, j_{1}' = j+2), (j_{2}=2, j_{2}' = k+2) \dots, (j_{n+1},j_{n+1}') \}$ of $S_{2(n+1)}$ as discussed above.  Let $\lambda = (\lambda_{1}, \dots, \lambda_{2n})$ with $\lambda_{1} \geq \lambda_{2} \geq \cdots \geq \lambda_{2n} \in \mathbb{Z}$.  Then we have the following term-evaluation:
\begin{multline*}
2^{n-1}(n-1)! \PV \int_{T} R_{\lambda, w}(x_{1}^{\pm 1}, \dots, x_{n-1}^{\pm 1}, \pm 1;t) \tilde \Delta_{K}^{(n-1)}(\pm t, \pm \sqrt{t}) \prod_{i=1}^{n-1} (1-\alpha x_{i}^{\pm 1}) dT  \\
= \frac{1+t}{2} \frac{\epsilon(w)}{2^{n-1}(1-t)^{n-1}} \prod_{1 \leq k \leq n+1} M_{j_{k},j_{k}'}^{\lambda}.
\end{multline*}
\end{claim}  
We prove the claim.  Let $\mu$ be such that $\lambda = \mu + \rho_{2n}$.  By Claim \ref{ame1} the above LHS is equal to 
\begin{multline*}
\begin{cases}  2^{n-1}(n-1)! (-1)^{\lambda_{j}+k-2} \frac{(1+t)}{2} \int R_{\tilde \lambda, \tilde w}(x_{1}^{\pm 1}, \dots, x_{n-1}^{\pm 1};t) \tilde \Delta_{K}^{(n-1)}(\pm 1, \pm \sqrt{t}) \prod_{i=1}^{n-1}(1- \alpha x_{i}^{\pm 1})dT& \text{$j<k$,} \\
2^{n-1}(n-1)! (-1)^{\lambda_{j}+k-1} \frac{(1+t)}{2} \int R_{\tilde \lambda, \tilde w}(x_{1}^{\pm 1}, \dots, x_{n-1}^{\pm 1};t) \tilde \Delta_{K}^{(n-1)}(\pm 1, \pm \sqrt{t}) \prod_{i=1}^{n-1}(1- \alpha x_{i}^{\pm 1})dT&\text{$j>k$,}
\end{cases}\\
= 2^{n-1} (n-1)! \frac{1+t}{2} (-1)^{j_{1}'+j_{2}'-1-c_{2}(w)} M_{1,j_{1}'}^{\lambda} M_{2, j_{2}'}^{\lambda} \\ \cdot \int_{T} R_{\tilde \lambda, \tilde w}(x_{1}^{\pm 1}, \dots, x_{n-1}^{\pm 1};t) \tilde \Delta_{K}^{(n-1)}(\pm 1, \pm \sqrt{t}) \prod_{i=1}^{n-1} (1-\alpha x_{i}^{\pm 1}) dT,
\end{multline*}
where $c_{2}(w)$ is $0$ if $j_{1}' > j_{2}'$ (i.e., $(1,j_{1}')$ and $(2, j_{2}')$ do not cross) and $1$ if they do.  Now we may use Claim \ref{apet} on the $(n-1)$-dimensional integral: let $\widehat{\mu}$ be the partition $\mu$ with parts $\mu_{j}$ and $\mu_{k}$ deleted; note that $\tilde \lambda$ and $\widehat{\mu} + \rho_{2n-2}$ have equivalent parts modulo $2$.  Using this, we find that the above is equal to 
\begin{multline*}
2^{n-1} (n-1)! \frac{1+t}{2} (-1)^{j_{1}'+j_{2}'-1-c_{2}(w)} M_{1,j_{1}'}^{\lambda} M_{2, j_{2}'}^{\lambda}  \\\cdot \int_{T} R_{\widehat{\mu} + \rho_{2n-2}, \tilde w}(x_{1}^{\pm 1}, \dots, x_{n-1}^{\pm 1};t) \tilde \Delta_{K}^{(n-1)}(\pm 1, \pm \sqrt{t}) \prod_{i=1}^{n-1} (1-\alpha x_{i}^{\pm 1}) dT \\
= \frac{1+t}{2} (-1)^{j_{1}'+j_{2}'-1-c_{2}(w)} M_{1,j_{1}'}^{\lambda}  M_{2, j_{2}'}^{\lambda} \frac{\epsilon(\tilde w)}{2^{n-1}(1-t)^{n-1}} \prod_{1 \leq k \leq n-1} a_{i_{k},i_{k}'}^{\widehat{\mu} + \rho_{2n-2}} \\
= \frac{1+t}{2} \frac{\epsilon(w)}{2^{n-1}(1-t)^{n-1}} \prod_{1 \leq k \leq n+1} M_{j_{k},j_{k}'}^{\lambda},
\end{multline*}
as desired.

Note that in particular this result shows that the integral of a matching is a term in $ \mathrm{Pf}[M]^{\lambda}(1+t)/2^{n}(1-t)^{n-1}$.

Now using the claim, we have
\begin{align*}
\int R_{\lambda}(x_{1}^{\pm 1}, \dots, x_{n-1}^{\pm 1},1,-1;t) \tilde \Delta_{K}^{(n-1)}(\pm t, \pm \sqrt{t}) \prod_{i=1}^{n-1} (1-\alpha x_{i}^{\pm 1}) dT 
= \frac{(1+t)}{2} \frac{1}{2^{n-1}(1-t)^{n-1}}\mathrm{Pf}[M]^{\lambda},
\end{align*}
since the terms of the Pfaffian are in bijection with the integrals of the pseudo-matchings. 

Finally, to prove the theorem, we use Proposition \ref{normalizations}(iv) to obtain
\begin{multline*}
\frac{(1-\alpha^{2})}{\int \tilde \Delta_{K}^{(n-1)}(\pm t, \pm \sqrt{t}) dT} \int P_{\lambda}(x_{1}^{\pm 1}, \dots, x_{n-1}^{\pm 1},1,-1;t) \tilde \Delta_{K}^{(n-1)}(\pm t, \pm \sqrt{t}) \prod_{i=1}^{n-1} (1-\alpha x_{i}^{\pm 1}) dT\\ 
= \frac{(1-\alpha^{2}) (1-t)(1-t^{2}) \cdots (1-t^{2n})}{v_{\lambda}(t)(1-t)^{n+1}} \frac{1}{2^{n}(1-t)^{n-1}}\mathrm{Pf}[M]^{\lambda}
= \frac{\phi_{2n}(t)}{v_{\lambda}(t)(1-t)^{2n}} \frac{(1-\alpha^{2})}{2^{n}}\mathrm{Pf}[M]^{\lambda}.
\end{multline*}

Following the computation in \cite[5.21]{FR} (but noting that they are missing a factor of $2$), $\mathrm{Pf}[M]^{\lambda}$ may be evaluated as
\begin{align*}
\frac{2^{n}}{(1-\alpha^{2})}\Big[ (-\alpha)^{\sum_{j=1}^{2n} [\lambda_{j} \text{ mod}2]} - (-\alpha)^{\sum_{j=1}^{2n} [(\lambda_{j}+1) \text{ mod}2]} \Big],
\end{align*}
which proves the theorem.  
\end{proof}

\begin{theorem}\label{apo}
Let $l(\lambda) \leq 2n+1$.  We have the following integral identity for $O^{+}(2n+1)$:
\begin{multline*}
 \frac{(1-\alpha)}{\int \tilde \Delta_{K}^{(n)}(t,-1,\pm \sqrt{t}) dT} \int P_{\lambda}(x_{1}^{\pm 1}, \dots, x_{n}^{\pm 1}, 1;t) \tilde \Delta_{K}^{(n)}(t,-1,\pm \sqrt{t})\prod_{i=1}^{n}(1-\alpha x_{i}^{\pm 1})dT \\
= \frac{\phi_{2n+1}(t)}{v_{\lambda}(t)(1-t)^{2n+1}} \Big[ (-\alpha)^{\# \text{ of odd parts of $\lambda$}} + (-\alpha)^{\# \text{ of even parts of $\lambda$}} \Big].
\end{multline*}
\end{theorem}

\begin{proof}
We use an argument analogous to the $O^{-}(2n)$ case.  We will first show the following:
\begin{align*}
\int R_{\lambda}(x_{1}^{\pm 1}, \dots, x_{n}^{\pm 1}, 1;t) \tilde \Delta_{K}^{(n)}(t,-1, \pm \sqrt{t}) \prod_{i=1}^{n} (1-\alpha x_{i}^{\pm 1}) dT 
= \frac{1}{2^{n}(1-t)^{n}} \mathrm{Pf}[M]^{\lambda},
\end{align*}
where the $2n+2 \times 2n+2$ antisymmetric matrix $[M]^{\lambda}$ is given by
\begin{align*}
\begin{cases} M_{1,k}^{\lambda} = 1 &\text{if $1 <k  \leq 2n+2$} \\
M_{j,k}^{\lambda} = a_{j-1,k-1}^{\lambda} &\text{if $2 \leq j \leq k \leq 2n+2$},
\end{cases}
\end{align*}
and as usual $[a_{j,k}]^{\lambda}$ is the $2n+1 \times 2n+1$ antisymmetric matrix specified by Theorem \ref{ape}.  The integral is a sum of $(2n+1)!$ terms, one for each permutation in $S_{2n+1}$.  But note that by symmetry we may restrict to pseudo-matchings in $S_{2n+1}$: those with $1$ anywhere but $x_{i}$ to the left of $x_{i}^{-1}$ for all $1 \leq i \leq n$, and $x_{i}$ to the left of $x_{j}$ for $1 \leq i<j \leq n$.  There are $(2n+1)!/2^{n}n!$ such pseudo-matchings, and for each there are exactly $2^{n}n!$ other permutations with identical integral value.

\begin{claim}\label{apo1}
Let $w$ be a fixed pseudo-matching with $1$ in position $k$, for some $1 \leq k \leq 2n+1$.  Then we have the following:
\begin{multline*}
 2^{n}n! \PV \int R_{\lambda, w}(x_{1}^{\pm 1}, \dots, x_{n}^{\pm 1},1;t) \tilde \Delta_{K}^{(n)}(t,-1,\pm \sqrt{t}) \prod_{i=1}^{n}(1- \alpha x_{i}^{\pm 1}) dT \\
= 2^{n} n! (-1)^{k-1} \PV \int R_{\tilde \lambda, \tilde w}(x_{1}^{\pm 1}, \dots, x_{n}^{\pm 1};t) \tilde \Delta_{K}^{(n)}(\pm 1, \pm \sqrt{t}) \prod_{i=1}^{n} (1-\alpha x_{i}^{\pm 1}) dT,
\end{multline*}
where $\tilde w$ is $w$ with $1$ deleted (in particular, a matching in $S_{2n}$) and $\tilde \lambda$ is $\lambda$ with $\lambda_{k}$ deleted and the parts to the left of $\lambda_{k}$ increased by $1$, i.e., 
\begin{align*}
\tilde \lambda &= (\lambda_{1} + 1, \dots, \lambda_{k-1}+1, \lambda_{k+1}, \dots, \lambda_{2n+1}).
\end{align*}
\end{claim}
We prove the claim; note that this proof is very similar to Claim \ref{ame1} for the $O^{-}(2n)$ case.  First, using (\ref{Kd}), we have
\begin{equation*}
2^{n}n! \tilde \Delta_{K}^{(n)}(t,-1, \pm \sqrt{t}) 
= \prod_{1 \leq i \leq n} \frac{1-x_{i}^{\pm 2}}{(1-tx_{i}^{\pm 1})(1+x_{i}^{\pm 1})(1-\sqrt{t}x_{i}^{\pm 1})(1+ \sqrt{t}x_{i}^{\pm 1})} \prod_{1 \leq i<j \leq n} \frac{1-x_{i}^{\pm 1}x_{j}^{\pm 1}}{1-tx_{i}^{\pm 1}x_{j}^{\pm 1}}.
\end{equation*}
Define the set $X = \{ (x_{i}^{\pm 1}, x_{j}^{\pm 1}): 1 \leq i \neq j \leq n \}$, and let $u_{\lambda, w}^{(n)}(x;t)$ be defined by
\begin{equation*}
R_{\lambda, w}(x_{1}^{\pm 1}, \dots, x_{n}^{\pm 1},1;t) = u_{\lambda, w}^{(n)}(x;t) \prod_{\substack{(z_{i},z_{j}) \in X: \\ z_{i} \prec_{w} z_{j}}} \frac{z_{i} - tz_{j}}{z_{i} - z_{j}}.
\end{equation*}
Also define $p_{1}$ and $\Delta_{1}$ by
\begin{equation*}
u_{\lambda, w}^{(n)}(x;t) \prod_{1 \leq i \leq n} \frac{1-x_{i}^{\pm 2}}{(1-tx_{i}^{\pm 1})(1+x_{i}^{\pm 1})(1-\sqrt{t}x_{i}^{\pm 1})(1+ \sqrt{t}x_{i}^{\pm 1})}  \prod_{i=1}^{n}(1- \alpha x_{i}^{\pm 1}) = p_{1} \prod_{i=1}^{n} \frac{1}{1-x_{i}^{-2}}
\end{equation*}
and 
\begin{equation*}
\prod_{1 \leq i<j \leq n} \frac{1-x_{i}^{\pm 1}x_{j}^{\pm 1}}{1-tx_{i}^{\pm 1}x_{j}^{\pm 1}}  \prod_{\substack{(z_{i},z_{j}) \in X: \\ z_{i} \prec_{w} z_{j}}} \frac{z_{i} - tz_{j}}{z_{i} - z_{j}} = \Delta_{1}. 
\end{equation*}
Note that
\begin{equation*}
R_{\lambda, w}(x_{1}^{\pm 1}, \dots, x_{n}^{\pm 1}, 1;t) \tilde \Delta_{K}^{(n)}(t,-1, \pm \sqrt{t})  \prod_{i=1}^{n}(1- \alpha x_{i}^{\pm 1}) = p_{1} \Delta_{1} \prod_{i=1}^{n} \frac{1}{1-x_{i}^{-2}}.
\end{equation*}
Define analogously $p_{2}$ and $\Delta_{2}$ using $R_{\tilde \lambda, \tilde w}(x_{1}^{\pm 1}, \dots, x_{n}^{\pm 1};t)$ and $\tilde \Delta_{K}^{(n)}(\pm 1, \pm \sqrt{t})$ instead of $R_{\lambda, w}^{(2n+1)}(x^{\pm 1},1;t)$ and $\tilde \Delta_{K}^{(n)}(t,-1, \pm \sqrt{t})$.  

Then note that $\Delta_{1} = \Delta_{2} := \Delta$.  Some computation shows that $\Delta(\pm 1, \dots, \pm 1, x_{i+1}, \dots, x_{n})$ is holomorphic in $x_{i+1}$ for all $0 \leq i \leq n-1$ and all $2^{i}$ combinations.  Further computations show that the function $p = p_{1} - (-1)^{k-1}p_{2}$ satisfies the conditions of Lemma \ref{pvl2}, so we have
\begin{align*}
\int p \cdot \Delta \cdot \prod_{i=1}^{n} \frac{1}{1-x_{i}^{-2}} dT &= 0
\end{align*}
or, 
\begin{align*}
\int p_{1} \cdot \Delta_{1} \cdot \prod_{i=1}^{n} \frac{1}{1-x_{i}^{-2}} dT &= (-1)^{k-1} \int p_{2} \cdot \Delta_{2} \cdot \prod_{i=1}^{n} \frac{1}{1-x_{i}^{-2}} dT,
\end{align*}
which proves the claim.

In keeping with the notation of the previous two theorems, we write $\{(k), (i_{1}, i_{1}'), \dots, (i_{n}, i_{n}') \}$ for the pseudo-matching $w$ with $1$ in position $k$ and $x_{k}$ in position $i_{k}$, $x_{k}^{-1}$ in position $i_{k}'$, for all $1 \leq k \leq n$.  We can extend this to a matching in $S_{2(n+1)}$ by $\{ (1, k+1), (i_{1} + 1, i_{1}'+1), \dots, (i_{n} + 1, i_{n}' + 1) \} = \{ (j_{1} = 1, j_{1}' = k+1), \dots, (j_{n+1}, j_{n+1}') \} $, with $i_{k} + 1 = j_{k+1}, i_{k'}+1 = j_{k+1}'$ for $1 \leq k \leq n$.  

\begin{claim}\label{apot}
Let $w = \{ (k), (i_{1}, i_{1}'), \dots, (i_{n}, i_{n}') \}$ be a pseudo-matching in $S_{2n+1}$, and extend it to a matching $\{ (j_{1} = 1, j_{1}' = k+1), \dots, (j_{n+1}, j_{n+1}') \}$ as discussed above.  Let $\lambda = (\lambda_{1}, \dots, \lambda_{2n+1})$ with $\lambda_{1} \geq \lambda_{2} \geq \cdots \geq \lambda_{2n+1} \in \mathbb{Z}$.  Then we have the following term-evaluation:
\begin{align*}
 2^{n}n! \PV \int_{T} R_{\lambda, w}(x_{1}^{\pm 1}, \dots, x_{n}^{\pm 1},1;t) \tilde \Delta_{K}^{(n)}(t,-1,\pm \sqrt{t}) \prod_{i=1}^{n} (1-\alpha x_{i}^{\pm 1}) dT 
= \frac{\epsilon(w)}{2^{n}(1-t)^{n}} \prod_{1 \leq k \leq n+1}  M_{j_{k}, j_{k}'}^{\lambda}.
\end{align*}
\end{claim}
We prove the claim.  Let $\mu$ be such that $\lambda = \mu + \rho_{2n+1}$.  By Claim \ref{apo1} the above LHS is equal to 
\begin{align*}
& 2^{n}n! (-1)^{k-1} \int_{T} R_{\tilde \lambda, \tilde w} (x_{1}^{\pm 1}, \dots, x_{n}^{\pm 1};t) \tilde \Delta_{K}^{(n)}(\pm 1, \pm \sqrt{t}) \prod_{i=1}^{n} (1-\alpha x_{i}^{\pm 1}) dT \\
&= 2^{n}n! (-1)^{j_{1}'-j_{1} + 1} M_{j_{1}, j_{1}'}^{\lambda} \int_{T} R_{\tilde \lambda, \tilde w} (x_{1}^{\pm 1}, \dots, x_{n}^{\pm 1};t) \tilde \Delta_{K}^{(n)}(\pm 1, \pm \sqrt{t}) \prod_{i=1}^{n} (1-\alpha x_{i}^{\pm 1}) dT. 
\end{align*}
Now we use Claim \ref{apet}: let $\widehat{\mu}$ be $\mu$ with part $\mu_{k}$ deleted; note $\tilde \lambda - 1^{2n} = \widehat{\mu} + \rho_{2n}$.  Using that result, the above is equal to 
\begin{multline*}
2^{n}n! (-1)^{j_{1}'-j_{1} + 1} M_{j_{1}, j_{1}'}^{\lambda} \int_{T} R_{\widehat{\mu} + \rho_{2n}, \tilde w} (x_{1}^{\pm 1}, \dots, x_{n}^{\pm 1};t) \tilde \Delta_{K}^{(n)}(\pm 1, \pm \sqrt{t}) \prod_{i=1}^{n} (1-\alpha x_{i}^{\pm 1}) dT \\
= (-1)^{j_{1}'-j_{1} + 1} M_{j_{1}, j_{1}'}^{\lambda} \frac{\epsilon(\tilde w)}{2^{n}(1-t)^{n}} \prod_{1 \leq k \leq n} a_{i_{k},i_{k}'}^{\widehat{\mu} + \rho_{2n}}= \frac{\epsilon(w)}{2^{n}(1-t)^{n}} \prod_{1 \leq k \leq n+1} M_{j_{k},j_{k}'}^{\lambda},
\end{multline*}
as desired.

Note that in particular this result shows that the integral of a matching is a term in $\frac{1}{2^{n}(1-t)^{n}} \mathrm{Pf}[M]^{\lambda}$.  

Now using the claim, we have
\begin{align*}
\int R_{\lambda}(x_{1}^{\pm 1}, \dots, x_{n}^{\pm 1}, 1;t) \tilde \Delta_{K}^{(n)}(t,-1, \pm \sqrt{t}) \prod_{i=1}^{n} (1-\alpha x_{i}^{\pm 1}) dT 
= \frac{1}{2^{n}(1-t)^{n}} \mathrm{Pf}[M]^{\lambda},
\end{align*}
since the terms of the Pfaffian are in bijection with the integrals of the pseudo-matchings.  

Finally, to prove the theorem, we use Proposition \ref{normalizations}(v) to obtain
\begin{multline*}
 \frac{(1-\alpha)}{\int \tilde \Delta_{K}^{(n)}(t,-1,\pm \sqrt{t}) dT} \int P_{\lambda}(x_{1}^{\pm 1}, \dots, x_{n}^{\pm 1}, 1;t) \tilde \Delta_{K}^{(n)}(t,-1,\pm \sqrt{t})\prod_{i=1}^{n}(1-\alpha x_{i}^{\pm 1})dT  \\
= \frac{(1-\alpha)\phi_{2n+1}(t)}{v_{\lambda}(t)(1-t)^{n+1}} \frac{1}{2^{n}(1-t)^{n}}\mathrm{Pf}[M]^{\lambda},
\end{multline*}
but by a change of basis $[M]^{\lambda}$ is equivalent to the one defined in \cite[5.24]{FR}, and that Pfaffian was computed to be
\begin{align*}
\frac{2^{n}}{(1-\alpha)}\Big[ (-\alpha)^{\sum_{j=1}^{2n+1} [\lambda_{j} \text{ mod}2]} + (-\alpha)^{\sum_{j=1}^{2n+1} [(\lambda_{j}+1) \text{ mod}2]} \Big],
\end{align*}
which proves the theorem.
\end{proof}

\begin{theorem}\label{amo}
Let $l(\lambda) \leq 2n+1$.  We have the following integral identity for $O^{-}(2n+1)$:
\begin{multline*}
 \frac{(1+\alpha)}{\int \tilde \Delta_{K}^{(n)}(1,-t,\pm \sqrt{t}) dT} \int P_{\lambda}(x_{1}^{\pm 1}, \dots, x_{n}^{\pm 1}, -1;t) \tilde \Delta_{K}^{(n)}(1,-t,\pm \sqrt{t})\prod_{i=1}^{n}(1-\alpha x_{i}^{\pm 1})dT \\
= \frac{\phi_{2n+1}(t)}{v_{\lambda}(t)(1-t)^{2n+1}} \Big[ (-\alpha)^{\# \text{ of odd parts of $\lambda$}} - (-\alpha)^{\# \text{ of even parts of $\lambda$}} \Big],
\end{multline*}
\end{theorem}
\begin{proof}
We obtain the $O^{-}(2n+1)$ integral from the $O^{+}(2n+1)$ integral.  See the discussion for the $O^{-}(2n+1)$ integral in the next section.  The upshot is that the $O^{-}(2n+1)$ integral is $(-1)^{|\lambda|}$ times the $O^{+}(2n+1)$ integral with parameter $-\alpha$.  Using Theorem \ref{apo}, we get
\begin{align*}
(-1)^{|\lambda|} \frac{\phi_{2n+1}(t)}{v_{\lambda}(t)(1-t)^{2n+1}} \Big[ \alpha^{\# \text{ of odd parts of $\lambda$}} + \alpha^{\# \text{ of even parts of $\lambda$}} \Big]. \\
\end{align*}
But note that $(-1)^{\lambda_{i}}$ is $-1$ if $\lambda_{i}$ is odd, and $1$ if $\lambda_{i}$ is even, so that $(-1)^{|\lambda|} = (-1)^{\# \text{ of odd parts of $\lambda$}}$.  Also,
\begin{align*}
 (-1)^{\# \text{ of odd parts of $\lambda$}} (-1)^{\# \text{ of even parts of $\lambda$}} = (-1)^{2n+1} 
 = -1.
\end{align*}
Combining these facts gives the result.  
\end{proof}

We briefly mention some existing results related to Theorems \ref{ape}, \ref{ame}, \ref{apo}, and \ref{amo}.  First, note that these four results are $t$-analogs of the results of Proposition $2$ of \cite{FR}.  For example, in the $O^{+}(2n)$ case, that result states
\begin{align*}
\langle \text{det}(1_{2n} + \alpha U) s_{\rho}(U)\rangle_{U \in O^{+}(2n)} = \frac{1}{2^{n-1}} \mathrm{Pf}[a_{jk}] 
= \alpha^{\sum_{j=1}^{2n} [\rho_{j} \text{ mod}2]} + \alpha^{\sum_{j=1}^{2n} [(\rho_{j}+1) \text{ mod}2]},
\end{align*}
where $\langle \cdot \rangle_{O^{+}(2n)}$ denotes the integral with respect to the eigenvalue density of the group $O^{+}(2n)$.  

Also, note that the $\alpha = 0$ case of these identities gives that the four integrals
\begin{align*}
\frac{1}{Z} \int P_{\lambda}(x_{1}^{\pm 1}, \dots, x_{n}^{\pm 1};t) \tilde \Delta_{K}^{(n)}(\pm 1, \pm \sqrt{t}) dT \\
\frac{1}{Z} \int P_{\lambda}(x_{1}^{\pm 1}, \dots, x_{n-1}^{\pm 1}, \pm 1;t) \tilde \Delta_{K}^{(n-1)}(\pm t, \pm \sqrt{t}) dT \\
\frac{1}{Z} \int P_{\lambda}(x_{1}^{\pm 1}, \dots, x_{n}^{\pm 1},1;t) \tilde \Delta_{K}^{(n)}(t,-1,\pm \sqrt{t}) dT \\
\frac{1}{Z} \int P_{\lambda}(x_{1}^{\pm 1}, \dots, x_{n}^{\pm 1}, -1;t) \tilde \Delta_{K}^{(n)}(1,-t, \pm \sqrt{t})dT 
\end{align*}
vanish unless all $2n$ or $2n+1$ (as appropriate) parts of $\lambda$ have the same parity (see Theorem 4.1 of \cite{RV}).  Here $Z$ is the normalization: it makes the integral equal to unity when $\lambda$ is the zero partition.

\section{$\alpha, \beta$ version}
In this section, we further generalize the identities of the previous section by using the Pieri rule to add an extra parameter $\beta$.  The values are given in terms of Rogers--Szeg\H{o} polynomials (\ref{RSpol}).
\begin{theorem}\label{absum}
We have the following integral identities:
\begin{enumerate}
\item for $O(2n)$
\begin{multline*}
\frac{1}{\int \tilde \Delta_{K}^{(n)}(\pm 1, \pm \sqrt{t})dT} \int P_{\mu}(x_{1}^{\pm 1}, \dots, x_{n}^{\pm 1};t) \tilde \Delta_{K}^{(n)}( \pm 1, \pm \sqrt{t}) \prod_{i=1}^{n} (1-\alpha x_{i}^{\pm 1}) (1-\beta x_{i}^{\pm 1}) dT\\
+ \frac{(1-\alpha^{2})(1-\beta^{2})}{\int \tilde \Delta_{K}^{(n-1)}(\pm t, \pm \sqrt{t})dT} \int P_{\mu}(x_{1}^{\pm 1}, \dots, x_{n-1}^{\pm 1},1,-1;t) \tilde \Delta_{K}^{(n-1)}(\pm t, \pm \sqrt{t}) \prod_{i=1}^{n-1} (1-\alpha x_{i}^{\pm 1})(1-\beta x_{i}^{\pm 1})dT\\
= \frac{2\phi_{2n}(t)}{v_{\mu}(t) (1-t)^{2n}}  \Big[ \Big(\prod_{i \geq 0} H_{m_{2i}(\mu)}(\alpha \beta;t) \prod_{i \geq 0} H_{m_{2i+1}(\mu)}(\beta/\alpha;t)\Big)(-\alpha)^{\# \text{ of odd parts of } \mu} \Big].
\end{multline*}
\item for $O(2n+1)$
\begin{multline*}
\frac{(1-\alpha)(1-\beta)}{\int \tilde \Delta_{K}^{(n)}(t,-1,\pm \sqrt{t}) dT } \int P_{\mu}(x_{1}^{\pm 1}, \dots, x_{n}^{\pm 1}, 1;t) \tilde \Delta_{K}^{(n)}(t,-1,\pm \sqrt{t})\prod_{i=1}^{n}(1-\alpha x_{i}^{\pm 1}) (1-\beta x_{i}^{\pm 1}) dT \\
+  \frac{(1+\alpha)(1+\beta)}{\int \tilde \Delta_{K}^{(n)}(1,-t, \pm \sqrt{t}) dT} \int P_{\mu}(x_{1}^{\pm 1}, \dots, x_{n}^{\pm 1}, -1;t) \tilde \Delta_{K}^{(n)}(1,-t,\pm \sqrt{t})\prod_{i=1}^{n}(1-\alpha x_{i}^{\pm 1})(1-\beta x_{i}^{\pm 1})dT \\
= \frac{2\phi_{2n+1}(t)}{v_{\mu}(t)(1-t)^{2n+1}}  \Big[ \Big(\prod_{i \geq 0} H_{m_{2i}(\mu)}(\alpha \beta;t) \prod_{i \geq 0} H_{m_{2i+1}(\mu)}(\beta/\alpha;t)\Big)(-\alpha)^{\# \text{ of odd parts of } \mu} \Big].
\end{multline*}
\end{enumerate}
\end{theorem}
\begin{proof}
The proof follows Warnaar's argument (Theorem 1.1 of \cite{W}), with the only difference being that we take into account zero parts in the computation whereas Warnaar's infinite version is concerned only with nonzero parts.  The basic method is to use the Pieri rule for $P_{\mu}(x;t)e_{r}(x)$ in combination with the results of the previous section (the sum of the results of Theorems \ref{ape}, \ref{ame} for $O(2n)$ and similarly Theorems \ref{apo}, \ref{amo} for $O(2n+1)$).  Note that Warnaar starts with the case $a=b=0$ in his notation (the orthogonal group case) and successively applies the Pieri rule two times, introducing a parameter each time.  Because we proved the $\alpha$ case in the previous section, we need only use the Pieri rule once.  
\end{proof}

\begin{theorem}\label{abindiv}
Write $\lambda = 0^{m_{0}(\lambda)} \; 1^{m_{1}(\lambda)} \; 2^{m_{2}(\lambda)} \cdots$, with total number of parts $2n$ or $2n+1$ as necessary.  Then we have the following integral identities for the components of the orthogonal group:
\begin{enumerate}
\item for $O^{+}(2n)$
\begin{multline*}
 \frac{1}{Z} \int P_{\lambda}(x_{1}^{\pm 1}, \dots, x_{n}^{\pm 1};t) \tilde \Delta_{K}^{(n)}( \pm 1, \pm \sqrt{t}) \prod_{i=1}^{n} (1-\alpha x_{i}^{\pm 1}) (1-\beta x_{i}^{\pm 1}) dT\\
= \frac{\phi_{2n}(t)}{v_{\lambda}(t) (1-t)^{2n}}  \Big[ \Big(\prod_{i \geq 0} H_{m_{2i}(\lambda)}(\alpha \beta;t) \prod_{i \geq 0} H_{m_{2i+1}(\lambda)}(\beta/\alpha;t)\Big)(-\alpha)^{\# \text{ of odd parts of } \lambda} \\
+  \Big(\prod_{i \geq 0} H_{m_{2i+1}(\lambda)}(\alpha \beta;t) \prod_{i \geq 0} H_{m_{2i}(\lambda)}(\beta/\alpha;t)\Big)(-\alpha)^{\# \text{ of even parts of } \lambda} \Big]\\
\end{multline*}
\item for $O^{-}(2n)$
\begin{multline*}
\frac{(1-\alpha^{2})(1-\beta^{2})}{Z} \int P_{\lambda}(x_{1}^{\pm 1}, \dots, x_{n-1}^{\pm 1},1,-1;t) \tilde \Delta_{K}^{(n-1)}(\pm t, \pm \sqrt{t}) \prod_{i=1}^{n-1} (1-\alpha x_{i}^{\pm 1})(1-\beta x_{i}^{\pm 1}) dT\\
= \frac{\phi_{2n}(t)}{v_{\lambda}(t) (1-t)^{2n}}  \Big[ \Big(\prod_{i \geq 0} H_{m_{2i}(\lambda)}(\alpha \beta;t) \prod_{i \geq 0} H_{m_{2i+1}(\lambda)}(\beta/\alpha;t)\Big)(-\alpha)^{\# \text{ of odd parts of } \lambda} \\
-  \Big(\prod_{i \geq 0} H_{m_{2i+1}(\lambda)}(\alpha \beta;t) \prod_{i \geq 0} H_{m_{2i}(\lambda)}(\beta/\alpha;t)\Big)(-\alpha)^{\# \text{ of even parts of } \lambda} \Big] \\
\end{multline*}
\item for $O^{+}(2n+1)$
\begin{multline*}
 \frac{(1-\alpha)(1-\beta)}{Z} \int P_{\lambda}(x_{1}^{\pm 1}, \dots, x_{n}^{\pm 1}, 1;t) \tilde \Delta_{K}^{(n)}(t,-1,\pm \sqrt{t})\prod_{i=1}^{n}(1-\alpha x_{i}^{\pm 1}) (1-\beta x_{i}^{\pm 1}) dT \\
= \frac{\phi_{2n+1}(t)}{v_{\lambda}(t)(1-t)^{2n+1}}  \Big[ \Big(\prod_{i \geq 0} H_{m_{2i}(\lambda)}(\alpha \beta;t) \prod_{i \geq 0 } H_{m_{2i+1}(\lambda)}(\beta/\alpha;t)\Big)(-\alpha)^{\# \text{ of odd parts of } \lambda} \\
+  \Big(\prod_{i \geq 0} H_{m_{2i+1}(\lambda)}(\alpha \beta;t) \prod_{i \geq 0} H_{m_{2i}(\lambda)}(\beta/\alpha;t)\Big)(-\alpha)^{\# \text{ of even parts of } \lambda} \Big]\\\
\end{multline*}
\item for $O^{-}(2n+1)$
\begin{multline*}
 \frac{(1+\alpha)(1+\beta)}{Z} \int P_{\lambda}(x_{1}^{\pm 1}, \dots, x_{n}^{\pm 1}, -1;t) \tilde \Delta_{K}^{(n)}(1,-t,\pm \sqrt{t})\prod_{i=1}^{n}(1-\alpha x_{i}^{\pm 1})(1-\beta x_{i}^{\pm 1})dT \\
= \frac{\phi_{2n+1}(t)}{v_{\lambda}(t)(1-t)^{2n+1}}  \Big[ \Big(\prod_{i \geq 0} H_{m_{2i}(\lambda)}(\alpha \beta;t) \prod_{i \geq 0} H_{m_{2i+1}(\lambda)}(\beta/\alpha;t)\Big)(-\alpha)^{\# \text{ of odd parts of } \lambda} \\
-  \Big(\prod_{i \geq 0} H_{m_{2i+1}(\lambda)}(\alpha \beta;t) \prod_{i \geq 0} H_{m_{2i}(\lambda)}(\beta/\alpha;t)\Big)(-\alpha)^{\# \text{ of even parts of } \lambda} \Big]\ ,
\end{multline*}
\end{enumerate}
where $Z$ is the normalization at $\alpha = 0, \beta=0$ and $\lambda = 0^{2n}, 0^{2n+1}$ as appropriate.
\end{theorem}

\begin{proof}
Note that the Hall--Littlewood polynomials satisfy the following property:
\begin{align*}
\Big( \prod_{i=1}^{l} z_{i} \Big) P_{\lambda}(z_{1}, \dots, z_{l};t) &= P_{\lambda + 1^{l}}(z_{1}, \dots, z_{l};t).\\
\end{align*}
So in the case $O(2n)$, for example, we have
\begin{align*}
P_{\mu}(x_{1}^{\pm 1}, \dots, x_{n}^{\pm 1};t) &= P_{\mu+1^{2n}}(x_{1}^{\pm 1}, \dots, x_{n}^{\pm 1};t) \\
P_{\mu}(x_{1}^{\pm 1}, \dots, x_{n-1}^{\pm 1},1,-1;t) &= -P_{\mu + 1^{2n}}(x_{1}^{\pm 1}, \dots, x_{n-1}^{\pm 1},1,-1;t). \\
\end{align*}
Thus,
\begin{multline*}
\frac{1}{\int \tilde \Delta_{K}^{(n)}(\pm 1, \pm \sqrt{t})dT} \int P_{\mu}(x_{1}^{\pm 1}, \dots, x_{n}^{\pm 1};t) \tilde \Delta_{K}^{(n)}( \pm 1, \pm \sqrt{t}) \prod_{i=1}^{n} (1-\alpha x_{i}^{\pm 1}) (1-\beta x_{i}^{\pm 1}) dT\\
- \frac{(1-\alpha^{2})(1-\beta^{2})}{\int \tilde \Delta_{K}^{(n-1)}(\pm t, \pm \sqrt{t})dT} \int P_{\mu}(x_{1}^{\pm 1}, \dots, x_{n-1}^{\pm 1},1,-1;t) \tilde \Delta_{K}^{(n-1)}(\pm t, \pm \sqrt{t}) \prod_{i=1}^{n-1} (1-\alpha x_{i}^{\pm 1})(1-\beta x_{i}^{\pm 1})dT\\
=\frac{1}{\int \tilde \Delta_{K}^{(n)}(\pm 1, \pm \sqrt{t})dT} \int P_{\mu + 1^{2n}}(x_{1}^{\pm 1}, \dots, x_{n}^{\pm 1};t) \tilde \Delta_{K}^{(n)}( \pm 1, \pm \sqrt{t}) \prod_{i=1}^{n} (1-\alpha x_{i}^{\pm 1}) (1-\beta x_{i}^{\pm 1}) dT\\
+ \frac{(1-\alpha^{2})(1-\beta^{2})}{\int \tilde \Delta_{K}^{(n-1)}(\pm t, \pm \sqrt{t})dT} \int P_{\mu + 1^{2n}}(x_{1}^{\pm 1}, \dots, x_{n-1}^{\pm 1},1,-1;t) \tilde \Delta_{K}^{(n-1)}(\pm t, \pm \sqrt{t}) \prod_{i=1}^{n-1} (1-\alpha x_{i}^{\pm 1})(1-\beta x_{i}^{\pm 1})dT\\
= \frac{2\phi_{2n}(t)}{v_{\mu + 1^{2n}}(t) (1-t)^{2n}}  \Big[ \Big(\prod_{i \geq 0} H_{m_{2i}(\mu + 1^{2n})}(\alpha \beta;t) \prod_{i \geq 0} H_{m_{2i+1}(\mu + 1^{2n})}(\beta/\alpha;t)\Big)(-\alpha)^{\# \text{ of odd parts of } \mu + 1^{2n}} \Big],
\end{multline*}
where the last equality follows from Theorem \ref{absum}(i).  Now note that $v_{\mu + 1^{2n}}(t) = v_{\mu}(t)$, $m_{i}(\mu + 1^{2n}) = m_{i-1}(\mu)$ for all $i \geq 1$, and the number of odd parts in $\mu +  1^{2n}$ is the same as the number of even parts in $\mu$.  Thus the above is equal to
\begin{equation*}
\frac{2\phi_{2n}(t)}{v_{\mu}(t) (1-t)^{2n}}  \Big[ \Big(\prod_{i \geq 0} H_{m_{2i+1}(\mu)}(\alpha \beta;t) \prod_{i \geq 0} H_{m_{2i}(\mu)}(\beta/\alpha;t)\Big)(-\alpha)^{\# \text{ of even parts of } \mu} \Big].
\end{equation*}
Then, taking the sum/difference of this equation and Theorem \ref{absum}(i), we obtain
\begin{multline*}
\frac{2}{\int \tilde \Delta_{K}^{(n)}(\pm 1, \pm \sqrt{t})dT} \int P_{\mu}(x_{1}^{\pm 1}, \dots, x_{n}^{\pm 1};t) \tilde \Delta_{K}^{(n)}( \pm 1, \pm \sqrt{t}) \prod_{i=1}^{n} (1-\alpha x_{i}^{\pm 1}) (1-\beta x_{i}^{\pm 1}) dT\\
= \frac{2\phi_{2n}(t)}{v_{\lambda}(t) (1-t)^{2n}}  \Big[ \Big(\prod_{i \geq 0} H_{m_{2i}(\lambda)}(\alpha \beta;t) \prod_{i \geq 0} H_{m_{2i+1}(\lambda)}(\beta/\alpha;t)\Big)(-\alpha)^{\# \text{ of odd parts of } \lambda} \\
+  \Big(\prod_{i \geq 0} H_{m_{2i+1}(\lambda)}(\alpha \beta;t) \prod_{i \geq 0} H_{m_{2i}(\lambda)}(\beta/\alpha;t)\Big)(-\alpha)^{\# \text{ of even parts of } \lambda} \Big],
\end{multline*}
and
\begin{multline*}
\frac{2(1-\alpha^{2})(1-\beta^{2})}{\int \tilde \Delta_{K}^{(n-1)}(\pm t, \pm \sqrt{t})} \int P_{\lambda}(x_{1}^{\pm 1}, \dots, x_{n-1}^{\pm 1},1,-1;t) \tilde \Delta_{K}^{(n-1)}(\pm t, \pm \sqrt{t}) \prod_{i=1}^{n-1} (1-\alpha x_{i}^{\pm 1})(1-\beta x_{i}^{\pm 1}) dT\\
= \frac{2\phi_{2n}(t)}{v_{\lambda}(t) (1-t)^{2n}}  \Big[ \Big(\prod_{i \geq 0} H_{m_{2i}(\lambda)}(\alpha \beta;t) \prod_{i \geq 0} H_{m_{2i+1}(\lambda)}(\beta/\alpha;t)\Big)(-\alpha)^{\# \text{ of odd parts of } \lambda} \\
-  \Big(\prod_{i \geq 0} H_{m_{2i+1}(\lambda)}(\alpha \beta;t) \prod_{i \geq 0} H_{m_{2i}(\lambda)}(\beta/\alpha;t)\Big)(-\alpha)^{\# \text{ of even parts of } \lambda} \Big], \\
\end{multline*}
as desired.  The $O(2n+1)$ result is analogous; use instead Theorem \ref{absum}(ii).  Note alternatively that as in the $\alpha$ case, we can obtain the $O^{-}(2n+1)$ integral directly from the $O^{+}(2n+1)$ integral, since the change of variables $x_{i} \rightarrow -x_{i}$ gives
\begin{multline*}
 \int P_{\lambda}(x_{1}^{\pm 1}, \dots, x_{n}^{\pm 1}, -1;t) \tilde \Delta_{K}^{(n)}(1,-t,\pm \sqrt{t}) \prod_{i=1}^{n} (1-\alpha x_{i}^{\pm 1})(1-\beta x_{i}^{\pm 1}) dT \\
= \int P_{\lambda}(-x_{1}^{\pm 1}, \dots, -x_{n}^{\pm 1}, -1;t)  \tilde \Delta_{K}^{(n)}(-1,t, \pm \sqrt{t}) \prod_{i=1}^{n} (1+ \alpha x_{i}^{\pm 1})(1+\beta x_{i}^{\pm 1}) dT\\
= (-1)^{|\lambda|}\int P_{\lambda}(x_{1}^{\pm 1}, \dots, x_{n}^{\pm 1}, 1;t) \tilde \Delta_{K}^{(n)}(-1, t, \pm \sqrt{t}) \prod_{i=1}^{n} (1+\alpha x_{i}^{\pm 1})(1+\beta x_{i}^{\pm 1}) dT,
\end{multline*}
and $\int \tilde \Delta_{K}^{(n)}(1,-t,\pm \sqrt{t}) dT = \int \tilde \Delta_{K}^{(n)}(-1,t,\pm \sqrt{t}) dT$, so that
\begin{multline*}
\frac{(1+\alpha)(1+\beta)}{\int \tilde \Delta_{K}^{(n)}(1,-t,\pm \sqrt{t}) dT}\int P_{\lambda}(x_{1}^{\pm 1}, \dots, x_{n}^{\pm 1}, -1;t) \tilde \Delta_{K}^{(n)}(1,-t,\pm \sqrt{t}) \prod_{i=1}^{n} (1-\alpha x_{i}^{\pm 1})(1-\beta x_{i}^{\pm 1}) dT \\
= \frac{(-1)^{|\lambda|}(1+\alpha)(1+\beta)}{\int \tilde \Delta_{K}^{(n)}(-1,t,\pm \sqrt{t})dT} \int P_{\lambda}(x_{1}^{\pm 1}, \dots, x_{n}^{\pm 1},1;t) \tilde \Delta_{K}^{(n)}(-1,t,\pm \sqrt{t}) \prod_{i=1}^{n}(1+\alpha x_{i}^{\pm 1})(1+\beta x_{i}^{\pm 1}) dT,
\end{multline*}
which is $(-1)^{|\lambda|}$ times the $O^{+}(2n+1)$ integral with parameters $-\alpha, -\beta$.  
\end{proof}
We remark that Theorem \ref{abindiv}(i) may be obtained using the direct method of the previous section.  One ultimately obtains a recursive formula, for which the Rogers--Szego polynomials are a solution.  However, this argument does not easily work for $O^{-}(2n), O^{+}(2n+1)$ and $O^{-}(2n+1)$.  Thus, it is more practical to use the Pieri rule to obtain the $O(l)$ ($l$ odd or even) integrals, and then solve for the components.

\section{Special Cases}
We will use the results of the previous section to prove some identities that correspond to particular values of $\alpha$ and $\beta$.  

\begin{corollary}\label{aminus1}
($\alpha = -1$)  We have the following identity:
\begin{align*}
\frac{1}{Z} \int P_{\lambda}^{(2n)}(x^{\pm 1};t) \tilde \Delta_{K}^{(n)}( \pm 1, \pm \sqrt{t}) \prod_{i=1}^{n}(1+x_{i}^{\pm 1})(1-\beta x_{i}^{\pm 1}) dT
= \frac{2\phi_{2n}(t)}{v_{\lambda}(t)(1-t)^{2n}}  \prod_{i \geq 0} H_{m_{i}(\lambda)}(-\beta;t),
\end{align*}
where the normalization $Z = \int \tilde \Delta_{K}^{(n)}(\pm 1, \pm \sqrt{t}) dT$.
\end{corollary}
\begin{proof}
Just put $\alpha = -1$ into \ref{abindiv}(i).  
\end{proof}

\begin{corollary}
($\alpha = -\beta$) We have the following identity:
\begin{multline*}
\frac{1}{Z} \int P_{\lambda}^{(2n)}(x^{\pm 1};t) \tilde \Delta_{K}^{(n)}(\pm 1, \pm \sqrt{t}) \prod_{i=1}^{n} (1-\alpha^{2} x_{i}^{\pm 2}) dT\\
= \frac{\phi_{2n}(t)}{v_{\lambda}(t)(1-t)^{2n}}\Big[ \Big(\prod_{i \geq 0} H_{m_{2i}(\lambda)}(-\alpha^{2};t) \prod_{i \geq 0}H_{m_{2i+1}(\lambda)}(-1;t)\Big) 
(-\alpha)^{\# \text{ of odd parts of } \lambda}
\\+  \Big(\prod_{i \geq 0} H_{m_{2i+1}(\lambda)}(-\alpha^{2};t) \prod_{i \geq 0} H_{m_{2i}(\lambda)}(-1;t)\Big)(-\alpha)^{\# \text{ of even parts of } \lambda} \Big],
\end{multline*}
where the normalization $Z = \int \tilde \Delta_{K}^{(n)}(\pm 1, \pm \sqrt{t}) dT$.  In particular, this vanishes unless all odd parts of $\lambda$ have even multiplicity, or all even parts of $\lambda$ have even multiplicity.  
\end{corollary}
\begin{proof}
Just put $\alpha = -\beta$ into Theorem \ref{abindiv}(i).  For the second part, we use \cite[1.10b]{W}: $H_{m}(-1;t)$ vanishes unless $m$ is even, in which case it is $(t;t^{2})_{m/2} = (1-t)(1-t^{3}) \cdots (1-t^{m-1})$.  
\end{proof}

\begin{corollary}
Symplectic Integral (see Theorem 4.1 of \cite{RV}).  We have the following identity:
\begin{align*}
\frac{1}{Z} \int P_{\lambda}(x_{1}^{\pm 1}, \dots, x_{n}^{\pm 1};t) \tilde \Delta_{K}^{(n)}(\pm \sqrt{t},0,0)dT = \frac{\phi_{n}(t^{2})}{(1-t^{2})^{n}v_{\mu}(t^{2})} =\frac{C^{0}_{\mu}(t^{2n};0,t^{2})}{C^{-}_{\mu}(t^{2};0,t^{2})},
\end{align*}
when $\lambda = \mu^{2}$ for some $\mu$ and $0$ otherwise (here the normalization $Z = \int \tilde \Delta_{K}^{(n)}( \pm \sqrt{t},0,0) dT$).  
\end{corollary}

\begin{proof}
Use the computation
\begin{align*}
\tilde \Delta_{K}^{(n)}(\pm \sqrt{t},0,0) &= \tilde \Delta_{K}^{(n)}(\pm 1, \pm \sqrt{t}) \prod_{1 \leq i \leq n} (1-\alpha x_{i}^{\pm 1})(1-\beta x_{i}^{\pm 1}) \big|_{\alpha = -1, \beta = 1},
\end{align*}
and Corollary \ref{aminus1} with $\beta = 1$.  The result then follows from \cite[1.10b]{W}: $H_{m_{i}(\lambda)}(-1;t)$ vanishes unless $m_{i}(\lambda)$ is even, in which case it is $(1-t)(1-t^{3}) \cdots (1-t^{m_{i}(\lambda)-1})$. 
\end{proof}

We remark that this integral identity may also be proved directly, using techniques similar to those used for the orthogonal group integrals of Section 4.  In fact, in this case, there are no poles on the unit circle so the analysis is much more straightforward.

\begin{corollary} \label{Kawanaka}
Kawanaka's identity (see \cite{Ka2}, \cite{Ka1}).  We have the following identity:
\begin{align*}
\frac{1}{Z} \int P_{\lambda}(x_{1}^{\pm 1}, \dots, x_{n}^{\pm 1};t) \tilde \Delta_{K}^{(n)}(1,\sqrt{t},0,0) &= \frac{\phi_{2n}(\sqrt{t})}{(1-\sqrt{t})^{2n}v_{\lambda}(\sqrt{t})} = \frac{C^0_\lambda(t^{n};0,\sqrt{t})}{C^-_\lambda(\sqrt{t};0,\sqrt{t})}
\end{align*}
(here the normalization $Z = \int \tilde \Delta_{K}^{(n)}(1,\sqrt{t},0,0) dT$).  
\end{corollary}
\begin{proof}
Use the computation
\begin{align*}
\tilde \Delta_{K}^{(n)}(1, \sqrt{t},0,0) &= \tilde \Delta_{K}^{(n)}(\pm 1, \pm \sqrt{t}) \prod_{1 \leq i \leq n} (1-\alpha x_{i}^{\pm 1})(1-\beta x_{i}^{\pm 1}) \big|_{\alpha = -1, \beta = -\sqrt{t}},
\end{align*}
and Corollary \ref{aminus1} with $\beta = -\sqrt{t}$.  The result then follows from \cite[1.10d]{W}: $H_{m}(\sqrt{t};t) = \prod_{j=1}^{m} (1+ (\sqrt{t})^{j})$.   
\end{proof}

\section{Limit $n \rightarrow \infty$}
In this section, we show that the $n \rightarrow \infty$ limit of Theorem \ref{abindiv}(i) in conjunction with the Cauchy identity gives Warnaar's identity (\cite[Theorem 1.1]{W}).  Thus, Theorem \ref{abindiv}(i) may be viewed as a finite dimensional analog of that particular generalized Littlewood identity.  
\begin{proposition}
(Gaussian result for $O^{+}(2n)$) For any symmetric function $f$,
\begin{align*}
\lim_{n \rightarrow \infty} \frac{\displaystyle \int f(x^{\pm 1}) \tilde \Delta_{K}^{(n)}(x;t;\pm 1, t_{2}, t_{3})dT}{\displaystyle \int \tilde \Delta_{K}^{(n)}(x;t;\pm 1, t_{2}, t_{3})dT} &= I_{G}(f; m; s),
\end{align*}
where $|t|, |t_{2}|, t_{3}| < 1$ and $m$ and $s$ are defined as follows:
\begin{align*}
m_{2k-1} &= \frac{t_{2}^{2k-1} + t_{3}^{2k-1}}{1-t^{2k-1}} \\
m_{2k} &= \frac{t_{2}^{2k} + t_{3}^{2k} + 1-t^{k}}{1-t^{2k}}\\
s_{k} &= \frac{k}{1-t^{k}}.
\end{align*}
Here $I_{G}(;m;s)$ is the Gaussian functional on symmetric functions defined by 
\begin{align*}
\int_{\mathbb{R}^{\deg(f)}} f \prod_{j=1}^{\deg(f)} (2 \pi s_{j})^{-1/2} e^{-(p_{j}-m_{j})^{2}/2s_{j}}dp_{j}.
\end{align*}
\end{proposition}
\begin{proof}
This is formally a special case of \cite[Theorem 7.17]{R}.  That proof relies on Theorem 6 of \cite{DS} and Section 8 of \cite{BR}.  The fact that two of the parameters $(t_{0}, \dots, t_{3})$ are $\pm 1$ makes that argument fail: however, replacing the symplectic group with $O^{+}(2n)$ resolves that issue.  
\end{proof}
Note that a similar argument would work for the components $O^{-}(2n), O^{+}(2n+1)$ and $O^{-}(2n+1)$.  

\begin{proposition}
We have the following:
\begin{multline*}
\lim_{n \rightarrow \infty} \frac{ \displaystyle \int \prod_{j,k} \frac{1-tx_{j}y_{k}^{\pm 1}}{1-x_{j}y_{k}^{\pm 1}} \displaystyle \prod_{k} (1-\alpha y_{k}^{\pm 1})(1-\beta y_{k}^{\pm 1}) \tilde \Delta_{K}^{(n)}(y;t; \pm 1, t_{2}, t_{3})dT}{\displaystyle \int \tilde \Delta_{K}^{(n)}(y;t; \pm 1, t_{2}, t_{3}) dT} \\
= \frac{(t_{2}\alpha, t_{3}\alpha, t_{2}\beta, t_{3}\beta;t)}{(\alpha^{2}t, \beta^{2}t; t^{2}) (\alpha\beta;t)} \prod_{j<k} \frac{1-tx_{j}x_{k}}{1-x_{j}x_{k}} \prod_{j} \frac{(1-tx_{j}^{2})(1-\alpha x_{j})(1-\beta x_{j})}{(1-t_{2}x_{j})(1-t_{3}x_{j})(1-x_{j})(1+x_{j})}.
\end{multline*}
\end{proposition}
\begin{proof}
Put
\begin{align*}
f = \prod_{j,k} \frac{1-tx_{j}y_{k}^{\pm 1}}{1-x_{j}y_{k}^{\pm 1}} \prod_{k} (1-\alpha y_{k}^{\pm 1})(1-\beta y_{k}^{\pm 1}) = \text{exp}\Big( \sum_{1 \leq k} \frac{p_{k}(x)p_{k}(y)(1-t^{k})}{k} - \frac{p_{k}(y)(\alpha^{k} + \beta^{k})}{k} \Big)
\end{align*}
(see \cite{Mac} for more details).  Then use the previous result, and complete the square in the Gaussian integral.  
\end{proof}

\begin{corollary}\label{cor}
We have the following identity in the limit:
\begin{multline*}
\lim_{n \rightarrow \infty} \frac{\displaystyle \int \prod_{j,k} \frac{1-tx_{j}y_{k}^{\pm 1}}{1-x_{j}y_{k}^{\pm 1}} \prod_{k} (1-\alpha y_{k}^{\pm 1})(1-\beta y_{k}^{\pm 1}) \tilde \Delta_{K}^{(n)}(y;t; \pm 1, \pm \sqrt{t})dT}{\displaystyle \int \tilde \Delta_{K}^{(n)}(y;t; \pm 1, \pm \sqrt{t}) dT} \\
= \frac{1}{(\alpha\beta;t)} \prod_{j<k} \frac{1-tx_{j}x_{k}}{1-x_{j}x_{k}} \prod_{j} \frac{(1-\alpha x_{j})(1-\beta x_{j})}{(1-x_{j})(1+x_{j})}.
\end{multline*}
\end{corollary}
\begin{proof}
Put $t_{2}, t_{3} = \pm \sqrt{t}$ in the previous result.  Also note that
\begin{equation*}
(\sqrt{t}\alpha;t)(-\sqrt{t}\alpha;t) = (t\alpha^{2};t^{2})
\end{equation*}
so that
\begin{align*}
\frac{(\sqrt{t}\alpha, -\sqrt{t}\alpha, \sqrt{t}\beta, -\sqrt{t}\beta;t)}{(\alpha^{2}t, \beta^{2}t;t^{2})} &= 1.
\end{align*}
\end{proof}

\begin{theorem}
We have the following formal identity (\cite{W} Theorem 1.1):
\begin{multline*}
\sum_{\lambda} P_{\lambda}(x;t)\Big[ \Big(\prod_{i > 0} H_{m_{2i}(\lambda)}(\alpha \beta;t) \prod_{i \geq 0} H_{m_{2i+1}(\lambda)}(\beta/\alpha;t)\Big)(-\alpha)^{\# \text{ of odd parts of } \lambda} \Big] \\
=  \prod_{j<k} \frac{1-tx_{j}x_{k}}{1-x_{j}x_{k}} \prod_{j} \frac{(1-\alpha x_{j})(1-\beta x_{j})}{(1-x_{j})(1+x_{j})}.
\end{multline*}
\end{theorem}

\begin{proof}
We prove the result for $|\alpha|, |\beta| < 1$, then use analytic continuation to obtain it for all $\alpha, \beta$.  We start with the Cauchy identity for Hall--Littlewood polynomials (\ref{Cauchyid}).  
Using this in the LHS of Corollary \ref{cor}, and multiplying both sides by $(\alpha \beta;t)$ gives
\begin{multline*}
(\alpha \beta;t) \sum_{\lambda} P_{\lambda}(x;t) \lim_{n \rightarrow \infty} \Big[ \frac{ b_{\lambda}(t)\int P_{\lambda}(y_{1}^{\pm 1}, \dots, y_{n}^{\pm 1};t) \prod_{k} (1-\alpha y_{k}^{\pm 1})(1-\beta y_{k}^{\pm 1}) \tilde \Delta_{K}^{(n)}(y;t; \pm 1, \pm \sqrt{t}) dT}{ \int \tilde \Delta_{K}^{(n)}(y;t;\pm 1, \pm \sqrt{t})dT}\Big] \\
= \prod_{j<k} \frac{1-tx_{j}x_{k}}{1-x_{j}x_{k}} \prod_{j} \frac{(1-\alpha x_{j})(1-\beta x_{j})}{(1-x_{j})(1+x_{j})}.
\end{multline*}
Now note that the quantity within the limit is the $\alpha, \beta$ version of the $O^{+}(2n)$ integral, see Theorem \ref{abindiv}(i).  Using that result, the above equation becomes
\begin{multline*}
(\alpha\beta;t)\sum_{\lambda}P_{\lambda}(x;t)\lim_{n \rightarrow \infty} \frac{b_{\lambda}(t)\phi_{2n}(t)}{v_{\lambda}(t) (1-t)^{2n}}  \Big[ \Big(\prod_{i \geq 0} H_{m_{2i}(\lambda)}(\alpha \beta;t) \prod_{i \geq 0} H_{m_{2i+1}(\lambda)}(\beta/\alpha;t)\Big)(-\alpha)^{\# \text{ of odd parts of } \lambda} \\
+  \Big(\prod_{i \geq 0} H_{m_{2i+1}(\lambda)}(\alpha \beta;t) \prod_{i \geq 0} H_{m_{2i}(\lambda)}(\beta/\alpha;t)\Big)(-\alpha)^{\# \text{ of even parts of } \lambda} \Big] \\
= \prod_{j<k} \frac{1-tx_{j}x_{k}}{1-x_{j}x_{k}} \prod_{j} \frac{(1-\alpha x_{j})(1-\beta x_{j})}{(1-x_{j})(1+x_{j})}.
\end{multline*}
But note that 
\begin{align*}
\frac{b_{\lambda}(t)}{v_{\lambda}(t)}&= \frac{(1-t)^{2n}}{\phi_{m_{0}(\lambda)}(t)},
\end{align*}
so that
\begin{align*}
\frac{b_{\lambda}(t)\phi_{2n}(t)}{v_{\lambda}(t)(1-t)^{2n}} &= \frac{\phi_{2n}(t)}{\phi_{m_{0}(\lambda)}(t)} = (1-t^{m_{0}(\lambda) +1}) \cdots (1-t^{2n}),
\end{align*}
which goes to $1$ as $m_{0}(\lambda), n \rightarrow \infty$.  Moreover, as $m_{0}(\lambda) \rightarrow \infty$, we have
\begin{multline*}
H_{m_{0}(\lambda)}(\alpha\beta;t) = \sum_{j=0}^{m_{0}(\lambda)} {\qbinom{m_{0}(\lambda)}{j}}_{t} (\alpha\beta)^{j} = \sum_{j=0}^{m_{0}(\lambda)} \frac{\phi_{m_{0}(\lambda)}(t)}{\phi_{j}(t)\phi_{m_{0}(\lambda)-j}(t)} (\alpha\beta)^{j}\\
= \sum_{j=0}^{m_{0}(\lambda)} \frac{(1-t^{m_{0}(\lambda)-j+1})(1-t^{m_{0}(\lambda)-j+2}) \cdots (1-t^{m_{0}(\lambda)})}{(1-t)(1-t^{2}) \cdots (1-t^{j})}(\alpha \beta)^{j} 
\rightarrow \sum_{j=0}^{\infty} \frac{(\alpha\beta)^{j}}{(t;t)_{j}}.
\end{multline*}
But for $|\alpha\beta|<1$, it is an identity that this is $1/(\alpha\beta;t)$.  

Finally, we show that the second term in the sum vanishes.  We must look at
\begin{align*}
\lim_{m_{0}(\lambda),k \rightarrow \infty }(-\alpha)^{k} H_{m_{0}(\lambda)}(\beta/\alpha;t) ,
\end{align*}
where $k$ is the number of even parts, so in particular $k \geq m_{0}(\lambda)$.  We have the following upper bound:
\begin{align*}
\lim_{m_{0}(\lambda) \rightarrow \infty} \alpha^{m_{0}(\lambda)} \sum_{j=0}^{m_{0}(\lambda)} \frac{(\beta/\alpha)^{j}}{(1-t)^{j}};
\end{align*}
the sum is geometric with ratio $\beta/\alpha(1-t)$.  Thus, this is equal to
\begin{equation*}
\lim_{m_{0}(\lambda) \rightarrow \infty} \alpha^{m_{0}(\lambda)} \frac{1-\Big(\frac{\beta}{\alpha(1-t)} \Big)^{m_{0}(\lambda)+1}}{1-\frac{\beta}{\alpha(1-t)} }
= \lim_{m_{0}(\lambda) \rightarrow \infty} \frac{\alpha^{m_{0}(\lambda)}- \frac{\beta^{m_{0}(\lambda)+1}}{\alpha (1-t)^{m_{0}(\lambda)+1}}}{1-\frac{\beta}{\alpha(1-t)} }.
\end{equation*}
But since $\alpha,\beta$ are sufficiently small (take $|\beta| < |1-t|$), this is zero, giving the result.  
\end{proof}

\section{Other Vanishing Results}

We introduce notation for dominant weights with negative parts: if $\mu, \nu$ are partitions with $l(\mu) + l(\nu) \leq n$ then $\mu \bar{\nu}$ is the dominant weight vector of $SL_{n} \times GL_{1}$, $\mu \bar{\nu} = (\mu_{1}, \dots, \mu_{l(\mu)},0, \cdots, 0, -\nu_{l(\nu)}, \dots, -\nu_{1})$.  Often, we will use $\lambda$ for a dominant weight with negative parts, i.e., $\lambda = \mu \bar{\nu}$.  

In this section, we prove four other vanishing identities from \cite{RV} and \cite{R}.  In all four cases, the structure of the partition that produces a nonvanishing integral is the same: opposite parts must add to zero ($\lambda_{i} + \lambda_{l+1-i} = 0$ for all $1 \leq i \leq l$, where $l$ is the total number of parts).  Note that an equivalent condition is that there exists a partition $\mu$ such that $\lambda = \mu \bar{\mu}$.    

We comment that the technique is similar to that of previous sections: we first use symmetries of the integrand to restrict to the term integrals associated to specific permutations.  Then, we obtain an inductive evaluation for the term integral, and use this to give a combinatorial formula for the total integral.  We mention that the first result corresponds to the symmetric space $(U(m+n), U(m) \times U(n))$ in the Schur case $t=0$.

\begin{theorem}(see \cite[Conjecture 3]{R})
Let $m$ and $n$ be integers with $0 \leq m \leq n$.  Then for a dominant weight $\lambda = \mu \bar{\nu}$ of $U(n+m)$,
\begin{align*}
\frac{1}{Z} \int_{T} P_{\mu \bar{\nu}}(x_{1}, \dots, x_{m}, y_{1}, \dots, y_{n};t) \frac{1}{n!m!} \prod_{1 \leq i \neq j \leq m} \frac{1-x_{i}x_{j}^{-1}}{1-tx_{i}x_{j}^{-1}} \prod_{1 \leq i \neq j \leq n} \frac{1-y_{i}y_{j}^{-1}}{1-ty_{i}y_{j}^{-1}} dT &= 0,
\end{align*}
unless $\mu = \nu$ and $l(\mu) \leq m$, in which case the integral is 
\begin{align*}
\frac{C^{0}_{\mu}(t^{n}, t^{m};0,t)}{C^{-}_{\mu}(t;0,t)C^{+}_{\mu}(t^{m+n-2}t;0,t)}.
\end{align*}
Here the normalization $Z$ is the integral for $\mu = \nu = 0$.
\end{theorem}
\begin{proof}
Note first that the integral is a sum of $(n+m)!$ terms, one for each element in $S_{n+m}$.  But by the symmetry of the integrand, we may restrict to the permutations with $x_{i}$ (resp. $y_{i}$) to the left of $x_{j}$ (resp. $y_{j}$) for $1 \leq i<j \leq m$ (resp. $1 \leq i<j \leq n$). Moreover, by symmetry we can deform the torus to
\begin{align*}
T= \{ |y| = 1 + \epsilon; |x| = 1 \},
\end{align*}
and preserve the integral.
Thus, we have 
\begin{multline*}
\int_{T} R_{\mu \bar{\nu}}(x^{(m)}, y^{(n)};t) \frac{1}{n!m!} \prod_{1 \leq i \neq j \leq m} \frac{1-x_{i}x_{j}^{-1}}{1-tx_{i}x_{j}^{-1}} \prod_{1 \leq i \neq j \leq n} \frac{1-y_{i}y_{j}^{-1}}{1-ty_{i}y_{j}^{-1}} dT \\
=  \sum_{\substack{w \in S_{n+m}\\ x_{i} \prec_{w} x_{j} \text{ for }1 \leq i<j \leq m \\ y_{i} \prec_{w} y_{j} \text{ for }1 \leq i<j \leq n}} \int_{T} R_{\mu \bar{\nu}, w}(x^{(m)}, y^{(n)};t)  \prod_{1 \leq i \neq j \leq m} \frac{1-x_{i}x_{j}^{-1}}{1-tx_{i}x_{j}^{-1}} \prod_{1 \leq i \neq j \leq n} \frac{1-y_{i}y_{j}^{-1}}{1-ty_{i}y_{j}^{-1}} dT
\end{multline*}
We first compute the normalization.
\begin{claim}\label{conj3norm}
We have
\begin{equation*}
 Z = \int_{T} P_{0^{n+m}} (x^{(m)}, y^{(n)};t) \frac{1}{n!m!}  \prod_{1 \leq i \neq j \leq m} \frac{1-x_{i}x_{j}^{-1}}{1-tx_{i}x_{j}^{-1}} \prod_{1 \leq i \neq j \leq n} \frac{1-y_{i}y_{j}^{-1}}{1-ty_{i}y_{j}^{-1}} dT= \frac{(1-t)^{m+n}}{\phi_{n}(t)\phi_{m}(t)}.
\end{equation*}
\end{claim}
Since
\begin{align*}
\frac{1}{v_{(0^{n+m})}(t)} &= \frac{(1-t)^{m+n}}{\phi_{m+n}(t)},
\end{align*}
this is equivalent to showing
\begin{align*}
&\int R_{0^{n+m}} (x^{(m)}, y^{(n)};t) \frac{1}{n!m!}  \prod_{1 \leq i \neq j \leq m} \frac{1-x_{i}x_{j}^{-1}}{1-tx_{i}x_{j}^{-1}} \prod_{1 \leq i \neq j \leq n} \frac{1-y_{i}y_{j}^{-1}}{1-ty_{i}y_{j}^{-1}} dT = \frac{\phi_{m+n}(t)}{\phi_{n}(t)\phi_{m}(t)}.
\end{align*}
We may use the above discussion to rewrite the LHS as a sum over suitable permutations.  Let $w \in S_{n+m}$ be a permutation with the $x$, $y$ variables in order and consider
\begin{equation*}
\int_{T} R_{0^{n+m}, w}(x^{(m)}, y^{(n)};t) \prod_{1 \leq i \neq j \leq m} \frac{1-x_{i}x_{j}^{-1}}{1-tx_{i}x_{j}^{-1}} \prod_{1 \leq i \neq j \leq n} \frac{1-y_{i}y_{j}^{-1}}{1-ty_{i}y_{j}^{-1}} dT.
\end{equation*}
Integrating with respect to $x_{1}, \dots, x_{m}, y_{1}, \dots, y_{n}$ in order shows that this is $t^{\# \text{inversions of } w}$, where inversions are in the sense of the multiset $M = \{ 0^{n}, 1^{m} \}$, and  we define $y_{1} \cdots y_{n} x_{1} \cdots x_{m}$ to have $0$ inversions.  But now by an identity of MacMahon
\begin{align*}
\sum_{\text{ multiset permutations $w$ of }\{0^{n}, 1^{m} \} } t^{\text{\# inversions of $w$}} &= {\qbinom{m+n}{n}}_{t} = \frac{\phi_{m+n}(t)}{\phi_{n}(t)\phi_{m}(t)},
\end{align*}
which proves the claim.  Note that we could also prove the claim by observing that
\begin{equation*}
\int_{T} P_{0^{n+m}} (x^{(m)}, y^{(n)};t) \frac{1}{n!m!}  \prod_{1 \leq i \neq j \leq m} \frac{1-x_{i}x_{j}^{-1}}{1-tx_{i}x_{j}^{-1}} \prod_{1 \leq i \neq j \leq n} \frac{1-y_{i}y_{j}^{-1}}{1-ty_{i}y_{j}^{-1}} dT = \frac{1}{n!m!} \int_{T} \tilde \Delta_{S}^{(m)}(x;t) \tilde \Delta_{S}^{(n)}(y;t) dT
\end{equation*}
and using the results of Theorem \ref{orthog}.

For convenience, from now on we will write
\begin{align*}
\Delta(x^{(m)}; y^{(n)};t) =  \prod_{1 \leq i \neq j \leq m} \frac{1-x_{i}x_{j}^{-1}}{1-tx_{i}x_{j}^{-1}} \prod_{1 \leq i \neq j \leq n} \frac{1-y_{i}y_{j}^{-1}}{1-ty_{i}y_{j}^{-1}} = \tilde \Delta_{S}^{(m)}(x;t) \tilde \Delta_{S}^{(n)}(y;t),
\end{align*}
for the density function.

\begin{claim}\label{conj3cl1}
Let $w \in S_{n+m}$ be a permutation of $\{x^{(m)},y^{(n)}\}$ with $x_{i} \prec_{w} x_{j}$ for all $1 \leq i<j \leq m$ and $y_{i} \prec_{w} y_{j}$ for all $1 \leq i<j \leq n$.  Suppose
\begin{align*}
\int_{T} R_{\mu\bar{\nu}, w}(x^{(m)},y^{(n)};t) \Delta(x^{(m)};y^{(n)};t) dT \neq 0.
\end{align*}
Then $w$ has $y_{1} \dots y_{l(\mu)}$ in first $l(\mu)$ positions, and $x_{m-l(\nu)+1} \dots x_{m}$ in the last $l(\nu)$ positions.  Consequently $l(\nu) \leq m$, $l(\mu) \leq n$.  
\end{claim}

We prove the claim.  We will first show that if, in $w(x,y)^{\mu \bar{\nu}}$, $x_{1}$ has exponent a strictly positive part, the integral is zero.  Indeed, one can compute that the integral restricted to the terms in $x_{1}$ is:
\begin{align*}
\int_{T_{1}} x_{1}^{\mu_{i}} \prod_{1 <i \leq m} \frac{x_{i}-x_{1}}{x_{i}-tx_{1}} \prod_{y_{j} \prec_{w} x_{1}} \frac{y_{j}-tx_{1}}{y_{j}-x_{1}} \prod_{x_{1} \prec_{w} y_{j}} \frac{x_{1}-ty_{j}}{x_{1}-y_{j}} dT &= 0,
\end{align*}
since by assumption $\mu_{i} > 0$.

Dually if in $w(x,y)^{\mu \bar{\nu}}$, $y_{n}$ has exponent a strictly negative part, we can show the integral is zero.  The integral restricted to the terms in in $y_{n}$ is:
\begin{multline*}
\int_{T_{1}} y_{n}^{\bar{\nu}_{i}} \prod_{1 \leq i < n} \frac{y_{n}-y_{i}}{y_{n}-ty_{i}} \prod_{x_{j} \prec_{w} y_{n}} \frac{x_{j}-ty_{n}}{x_{j}-y_{n}} \prod_{y_{n} \prec_{w} x_{j}} \frac{y_{n}-tx_{j}}{y_{n}-x_{j}}dT \\
= \int_{T: |x| > |y|}  y_{n}^{-\bar{\nu}_{i}} \prod_{1 \leq i < n} \frac{y_{i}-y_{n}}{y_{i}-ty_{n}} \prod_{x_{j} \prec_{w} y_{n}} \frac{y_{n}-tx_{j}}{y_{n}-x_{j}} \prod_{y_{n} \prec_{w} x_{j}} \frac{x_{j}-ty_{n}}{x_{j}-y_{n}}dT,
\end{multline*}
where in the second step we have inverted all variables which preserves the integral.  But now by assumption $\bar{\nu}_{i} < 0$, so integrating with respect to $y_{n}$ gives that the above integral is zero.  This gives the desired structure of $w$ to have nonvanishing associated integral.

\begin{claim}\label{conj3cl2}
Let $w \in S_{n+m}$ be a permutation of $\{x^{(m)},y^{(n)}\}$ with $x_{i} \prec_{w} x_{j}$ for all $1 \leq i<j \leq m$ and $y_{i} \prec_{w} y_{j}$ for all $1 \leq i<j \leq n$.  Suppose also that $y_{1}, \dots, y_{l(\mu)}$ are in the first $l(\mu)$ positions and $x_{m-l(\nu)+1}, \dots, x_{m}$ are in the last $l(\nu)$ positions.  

Let $l(\mu) > 0$.  Then we have the following formula for the term integral associated to $w$:
\begin{multline*}
\int_{T} R_{\mu \bar{\nu}, w}(x^{(m)},y^{(n)};t) \Delta(x^{(m)};y^{(n)};t) dT \\
= (1-t) \Big(\displaystyle \sum_{\substack{i: \\ \lambda_{1} + \lambda_{i} = 0}} t^{n+m - i} \Big) \int R_{\widehat{\lambda}, \widehat{w}}(x^{(m-1)}, y^{(n-1)};t) \Delta(x^{(m-1)};y^{(n-1)};t) dT  
\end{multline*}
where $\widehat{w}$ is $w$ with $y_{1}, x_{m}$ deleted and $\widehat{\lambda}$ is $\lambda$ with $\lambda_{1}$ and $\lambda_{i}$ deleted (where index $i$ is such that $\lambda_{1} + \lambda_{i} = 0$).

Similarly, if $l(\nu) > 0$, we have 
\begin{multline*}
\int_{T} R_{\mu \bar{\nu}, w}(x^{(m)},y^{(n)};t) \Delta(x^{(m)};y^{(n)};t) dT \\
= (1-t) \Big(\displaystyle \sum_{\substack{i: \\  \lambda_{i} + \lambda_{n+m} = 0}} t^{i-1} \Big) \int R_{\widehat{\lambda}, \widehat{w}}(x^{(m-1)}, y^{(n-1)};t) \Delta(x^{(m-1)};y^{(n-1)};t) dT  
\end{multline*}
where $\widehat{w}$ is $w$ with $y_{1}, x_{m}$ deleted and $\widehat{\lambda}$ is $\lambda$ with $\lambda_{i}$ and $\lambda_{n+m}$ deleted (where index $i$ is such that $\lambda_{i} + \lambda_{n+m} = 0$).

\end{claim}

For the first statement, integrate with respect to $y_{1}$.  We have the following integral restricted to the terms involving $y_{1}$: 
\begin{align*}
\int_{T_{1}} y_{1}^{\lambda_{1}} \prod_{1<i \leq n} \frac{y_{i}-y_{1}}{y_{i}-ty_{1}} \prod_{1 \leq j \leq m} \frac{y_{1}-tx_{j}}{y_{1} - x_{j}} dT,
\end{align*}
with $\lambda_{1} = \mu_{1} > 0$.  Evaluating gives a sum of $m$ terms, one for each residue $y_{1} = x_{j}$.  We consider one of these residues: suppose $x_{j}$ is in position $i$, then the resulting integral in $x_{j}$ is:
\begin{multline*}
(1-t) \int_{T_{1}} \negthickspace\negthickspace x_{j}^{\lambda_{1}+\lambda_{i}}  \prod_{1<i \leq n} \frac{y_{i} - x_{j}}{y_{i}-tx_{j}} \prod_{i \neq j} \frac{x_{j}-tx_{i}}{x_{j}-x_{i}} \prod_{\substack{y_{i} \prec_{w} x_{j} \\ y_{i} \neq y_{1}}} \frac{y_{i}-tx_{j}}{y_{i}-x_{j}} \prod_{x_{j} \prec_{w} y_{i}} \frac{x_{j}-ty_{i}}{x_{j}-y_{i}} \prod_{i<j} \frac{x_{j}-x_{i}}{x_{j}-tx_{i}} \prod_{j<i} \frac{x_{i}-x_{j}}{x_{i}-tx_{j}} dT \\
= (1-t) \int_{T_{1}} x_{j}^{\lambda_{1} + \lambda_{i}} \prod_{x_{j} \prec_{w} y_{i}}(-1) \frac{x_{j}-ty_{i}}{y_{i}-tx_{j}} \prod_{j<i} (-1) \frac{x_{j}-tx_{i}}{x_{i}-tx_{j}}dT,
\end{multline*}
where we may assume $\lambda_{i} \leq 0$, by the structure of $w$.  Note first that if $\lambda_{1} + \lambda_{i} >0$, the integral is zero.  One can similarly argue that the term integral is zero if $\lambda_{1} + \lambda_{i} <0$ (use $\lambda_{n+m} + \lambda_{k} < 0$ for any $1 \leq k <n+m$ and integrate with respect to $x_{m}$, and take the residue at any $x_{m} = y_{i}$).  Thus for a nonvanishing residue term we must have $\lambda_{1} = -\lambda_{i}$, and in this case one can verify that the above integral evaluates to
\begin{align*}
(1-t) t^{| \{ z: x_{j} \prec_{w} z \}| } = (1-t) t^{n+m-i},
\end{align*}
as desired.

The second statement is analogous, except integrate with respect to $x_{m}$ instead of $y_{1}$, and invert all variables.  This proves the claim.

Thus,
\begin{align*}
\int_{T} R_{\mu \bar{\nu}, w}(x^{(m)}, y^{(n)};t) \Delta(x^{(m)};y^{(n)};t) dT &= 0
\end{align*}
unless $\mu = \nu$ and $l(\mu) \leq m$, which gives the vanishing part of the theorem.  For the second part, suppose $\mu = \nu$ and $l(\mu) \leq m$.  Then by the above claims,
\begin{multline*}
\int_{T} R_{\mu\bar{\mu}, w}(x^{(m)}, y^{(n)};t) \Delta(x^{(m)};y^{(n)};t) dT \\
= (1-t)^{l(\mu)} v_{\mu+}(t) \int R_{0^{(n-l(\mu)) + (m-l(\nu))}, \delta}(x^{(m-l(\nu))}, y^{(n-l(\mu))};t)\Delta(x^{(m-l(\nu))};y^{(n-l(\mu))};t) dT  
\end{multline*}
if $w = y_{1} \dots y_{l(\mu)} \delta x_{m-l(\nu)+1} \dots x_{m}$ for some permutation $\delta$ of $\{y_{l(\mu)+1}, \dots, y_{n}, x_{1}, \dots, x_{m-l(\nu)} \}$, and $0$ otherwise.

By Claim \ref{conj3norm}, we have  
\begin{align*}
 \int R_{0^{(n-l(\mu)) + (m-l(\mu))}}(x^{(m-l(\mu))}, y^{(n-l(\mu))};t)\frac{\Delta(x^{(m-l(\mu))};y^{(n-l(\mu))};t)}{(m-l(\mu))!(n-l(\mu))!} dT  = \qbinom{m+n-2l(\mu)}{n-l(\mu)}_{t}.
\end{align*}

So we have
\begin{align*}
\int_{T} P_{\mu \bar{\mu}}(x^{(m)}, y^{(n)};t) \frac{1}{n!m!} \Delta(x^{(m)};y^{(n)};t) dT &= \frac{1}{v_{\mu\bar{\mu}}(t)} (1-t)^{l(\mu)} v_{\mu+}(t) \qbinom{m+n-2l(\mu)}{n-l(\mu)}_{t}.
\end{align*}

Noting that $v_{\mu\bar{\mu}}(t) = v_{\mu+}(t)^{2} v_{(0^{m+n-2l(\mu)})}(t)$ and multiplying by the reciprocal of the normalization gives
\begin{multline*}
\frac{1}{Z}\int_{T} P_{\mu \bar{\mu}}(x^{(m)}, y^{(n)};t)\frac{1}{n!m!} \Delta(x^{(m)};y^{(n)};t) dT= \frac{\phi_{n}(t)\phi_{m}(t)}{(1-t)^{m+n}}\frac{(1-t)^{l(\mu)}}{v_{\mu+}(t)v_{(0^{m+n-2l(\mu)})}(t)} \qbinom{m+n-2l(\mu)}{n-l(\mu)}_{t} \\
= (1-t^{n-l(\mu)+1}) \cdots (1-t^{n}) (1-t^{m-l(\mu)+1}) \cdots (1-t^{m}) \frac{ \phi_{m+n-2l(\mu)}(t) }{(1-t)^{m+n-l(\mu)}v_{\mu+}(t) v_{(0^{m+n-2l(\mu)})}(t)} \\
=\frac{(1-t^{n-l(\mu)+1}) (1-t^{n-l(\mu)+2}) \cdots (1-t^{n}) (1-t^{m-l(\mu)+1})(1-t^{m-l(\mu)+2}) \cdots (1-t^{m})}{(1-t)^{l(\mu)} v_{\mu+}(t)},
\end{multline*}
where the last equality follows from the definition of $v_{(0^{m+n-2l(\mu)})}$.
One can check from the definition of the $C$-symbols that 
\begin{align*}
C_{\mu}^{+}(t^{m+n-2}t;0,t) &= 1 \\
C^{-}_{\mu}(t;0,t) &= v_{\mu+}(t) (1-t)^{l(\mu)} \\
C^{0}_{\mu}(t^{n}, t^{m} ;0,t) &= \prod_{1 \leq i \leq l(\mu)} (1-t^{n+1-i})(1-t^{m+1-i}),
\end{align*}
so that our formula gives
\begin{align*}
\frac{C^{0}_{\mu}(t^{n},t^{m};0,t)}{C^{-}_{\mu}(t;0,t)C^{+}_{\mu}(t^{m+n-2}t;0,t)},
\end{align*}
as desired.
\end{proof}

\begin{theorem}
(see \cite[Conjecture 5]{R}) Let $n \geq 0$ be an integer and $\lambda = \mu \bar{\nu}$ a dominant weight of $U(2n)$.  Then 
\begin{align*}
\frac{1}{Z} \int_{T} P_{\mu \bar{\nu}}(x_{1}, \dots, x_{n}, y_{1}, \dots, y_{n};t) \frac{1}{(n!)^{2}} \prod_{1 \leq i,j \leq n} \frac{1}{(1-tx_{i}y_{j}^{-1})(1-ty_{i}x_{j}^{-1})} \prod_{1 \leq i \neq j \leq n} (1-x_{i}x_{j}^{-1})(1-y_{i}y_{j}^{-1})dT,
\end{align*}
is equal to 0 unless $\mu = \nu$, in which case the integral is
\begin{align*}
\frac{C^{0}_{\mu}(t^{n}, -t^{n};0,t)}{C^{-}_{\mu}(t;0,t)C^{+}_{\mu}(t^{2n-2}t;0,t)}.
\end{align*}
Here the normalization $Z$ is the integral for $\mu = \nu = 0$.
\end{theorem}
\begin{proof}
Note first that the integral is a sum of $(2n)!$ terms, one for each element in $S_{2n}$.  But by the symmetry of the integrand, we may restrict to the permutations with $x_{i}$ (resp. $y_{i}$) to the left of $x_{j}$ (resp. $y_{j}$) for $1 \leq i<j \leq n$.  By symmetry, we can deform the torus to 
\begin{align*}
T = \{ |y| = 1 + \epsilon; |x| = 1 \}.
\end{align*}
For convenience, we will write $\Delta(x^{(n)};y^{(n)};t)$ for the density
\begin{align*}
\prod_{1 \leq i,j \leq n} \frac{1}{(1-tx_{i}y_{j}^{-1})(1-ty_{i}x_{j}^{-1})} \prod_{1 \leq i \neq j \leq n} (1-x_{i}x_{j}^{-1})(1-y_{i}y_{j}^{-1}).
\end{align*}  
We first compute the normalization.
\begin{claim}\label{conj5norm}
We have
\begin{align*}
 Z = \int_{T} P_{0^{2n}}(x^{(n)}, y^{(n)};t) \frac{1}{(n!)^{2}} \Delta(x^{(n)};y^{(n)};t) dT = \frac{1}{\phi_{n}(t^{2})}.
\end{align*}
\end{claim}
By the definition of $v_{(0^{2n})}(t)$, this is equivalent to showing 
\begin{align*}
&\int_{T} R_{0^{2n}}(x^{(n)}, y^{(n)};t) \frac{1}{(n!)^{2}} \Delta(x^{(n)};y^{(n)};t) dT 
 = \frac{\phi_{2n}(t)}{(1-t)^{2n}\phi_{n}(t^{2})}.
\end{align*}
We prove this statement by induction on $n$.  For $n=1$, we have 
\begin{align*}
\int_{T}  \frac{x_{1}y_{1}}{(x_{1}-y_{1})(y_{1}-tx_{1})}dT = 0
\end{align*}
and
\begin{align*}
\int_{T} \frac{x_{1}y_{1}}{(y_{1}-x_{1})(x_{1}-ty_{1})} dT = \frac{1}{1-t} 
= \frac{\phi_{2}(t)}{(1-t)^{2}\phi_{1}(t^{2})}
\end{align*}
as desired.  Now suppose the claim holds for $n-1$; with this assumption we show that it holds for $n$.

Consider permutations $w$ with $x_{1}$ first.  We claim $\int_{T} R_{\mu \bar{\nu}, w}(x^{(n)},y^{(n)};t) \Delta(x^{(n)};y^{(n)};t) dT =0$.  Indeed, we have the following integral restricting to the terms in $x_{1}$:
\begin{multline*}
 \int_{T_{1}} \prod_{1 \leq i \leq n} \frac{x_{1} - ty_{i}}{x_{1}-y_{i}} \prod_{1<i \leq n} \frac{x_{1}-tx_{i}}{x_{1}-x_{i}} \prod_{1 \leq j \leq n} \frac{x_{1}y_{j}}{(y_{j}-tx_{1})(x_{1}-ty_{j})} \prod_{1<j \leq n} \frac{(x_{j}-x_{1})(x_{1}-x_{j})}{x_{1}x_{j}}dT \\
= \int_{T_{1}} \prod_{1 \leq j \leq n} \frac{x_{1}y_{j}}{(x_{1}-y_{j})(y_{j}-tx_{1})} \prod_{1<j \leq n} \frac{(x_{1}-tx_{j})(x_{j}-x_{1})}{x_{1}x_{j}} dT \\
= \int_{T_{1}} x_{1} \prod_{1 \leq j \leq n} \frac{1}{(x_{1}-y_{j})(y_{j}-tx_{1})} \prod_{1<j \leq n} (x_{1}-tx_{j})(x_{j}-x_{1}) dT 
= 0.
\end{multline*}
Thus, we may suppose $y_{1}$ occurs first in $w$.  A similar calculation for the integral restricting to terms in $y_{1}$ yields: 
\begin{align*}
\int_{T_{1}} y_{1} \prod_{1<j \leq n} (y_{1} - ty_{j})(y_{j}-y_{1}) \prod_{1 \leq i \leq n} \frac{1}{(y_{1} - x_{i})(x_{i}-ty_{1})}dT.
\end{align*}
We may evaluate this as the sum of $n$ residues, one for each $y_{1} = x_{i}$ for $1 \leq i \leq n$.  We compute the residue at $y_{1} = x_{i}$, and look at the resulting integral in $x_{i}$:
\begin{multline*}
\frac{1}{1-t}\int_{T_{1}}  \prod_{1<j \leq n} (x_{i}-ty_{j})(y_{j}-x_{i}) \prod_{j \neq i} \frac{1}{(x_{i}-x_{j})(x_{j}-tx_{i})} 
 \prod_{i' < i} (x_{i'} - tx_{i})(x_{i}-x_{i'}) \prod_{i<i''} (x_{i}-tx_{i''})(x_{i''}-x_{i}) \\
 \cdot \prod_{x_{i} \prec_{w} y_{j}} \frac{1}{(x_{i}-y_{j})(y_{j}-tx_{i})} \prod_{\substack{y_{j} \prec_{w} x_{i} \\ y_j \neq y_{1}}} \frac{1}{(y_{j}-x_{i})(x_{i}-ty_{j})} dT 
= \frac{1}{1-t}\int_{T_{1}}  \prod_{i< i''}  \frac{(tx_{i''}-x_i)}{(x_{i''}-tx_{i})} \prod_{x_{i} \prec_{w} y_{j}}  \frac{(ty_{j}-x_i)}{(y_{j}-tx_{i})} dT.
\end{multline*}
But, letting $2 \leq k \leq 2n$ be the position of $x_{i}$ in $w$, this evaluates to 
\begin{align*}
\frac{1}{1-t} \prod_{i< i''} t \prod_{x_{i} \prec_{w} y_{j}} t = \frac{t^{2n-k}}{1-t}.
\end{align*}
Thus, varying over all such permutations with $y_{1}$ first gives a factor of 
\begin{align*}
\frac{1}{1-t} (t^{2n-2} + t^{2n-3} + \cdots +t+1) &= \frac{(1-t^{2n-1})}{(1-t)^{2}}.
\end{align*}
Note that permutations of $\{y_{1}, \dots, y_{n}, x_{1}, \dots, x_{n}\}$ with $y_{1}$ in position $1$ and $x_{i}$ in position $k$ are in bijection with permutations of $\{y_{2}, \dots, y_{n}, x_{1}, \dots, \widehat{x_{i}}, \dots, x_{n}\}$.  So using the induction hypothesis, the total integral evaluates to
\begin{align*}
\frac{(1-t^{2n-1})}{(1-t)^{2}} \frac{\phi_{2(n-1)}(t)}{(1-t)^{2(n-1)}\phi_{n-1}(t^{2})} &= \frac{\phi_{2n}(t)}{(1-t)^{2n} \phi_{n}(t^{2})},
\end{align*}
as desired.  

Note that the density is not of a standard form (i.e., as a product of Koornwinder or Selberg densities), so we cannot appeal to an earlier result (compare with Claim \ref{conj3norm}).

\begin{claim}
Let $w \in S_{2n}$ a permutation of $\{x^{(n)},y^{(n)}\}$ with $x_{i} \prec_{w} x_{j}$ for all $1 \leq i<j \leq n$ and $y_{i} \prec_{w} y_{j}$ for all $1 \leq i<j \leq n$.  Suppose
\begin{align*}
\int_{T} R_{\mu \bar{\nu},w}(x^{(n)},y^{(n)};t) \Delta(x^{(n)};y^{(n)};t) dT \neq 0.
\end{align*}
Then $w$ has $y_{1} \dots y_{l(\mu)}$ in the first $l(\mu)$ coordinates, and $x_{n-l(\nu)+1} \dots x_{n}$ in the last $l(\nu)$ coordinates.  Consequently $l(\nu) \leq n, l(\mu) \leq n$.  
\end{claim}
The proof is analogous to Claim \ref{conj3cl1} of the previous theorem.

\begin{claim}
Let $w \in S_{2n}$ be a permutation of $\{x^{(n)},y^{(n)}\}$ with $x_{i} \prec_{w} x_{j}$ for all $1 \leq i<j \leq n$ and $y_{i} \prec_{w} y_{j}$ for all $1 \leq i<j \leq n$.  Suppose also that $y_{1}, \dots, y_{l(\mu)}$ are in the first $l(\mu)$ coordinates, and $x_{n-l(\nu)+1} \dots x_{n}$ in the last $l(\nu)$ coordinates.  

Let $l(\mu) > 0$.  Then we have the following formula for the term integral associated to $w$:
\begin{multline*}
\int_{T} R_{\mu \bar{\nu}, w}(x^{(n)},y^{(n)};t) \Delta(x^{(n)};y^{(n)};t) dT \\
= \frac{1}{1-t} \Big(\displaystyle \sum_{\substack{i : \\ \lambda_{1} + \lambda_{i} = 0}} t^{2n-i} \Big) \int R_{\widehat{\lambda}, \widehat{w}}(x^{(n-1)}, y^{(n-1)};t) \Delta(x^{(n-1)};y^{(n-1)};t) dT 
 \end{multline*}
where $\widehat{w}$ is $w$ with $y_{1}, x_{n}$ deleted and $\widehat{\lambda}$ is $\lambda$ with $\lambda_{1}$ and $\lambda_{i}$ deleted (where the index $i$ is such that $\lambda_{1} + \lambda_{i} = 0$).

Similarly, if $l(\nu) > 0$, we have 
\begin{multline*}
\int_{T} R_{\mu \bar{\nu}, w}(x^{(n)},y^{(n)};t) \Delta(x^{(n)};y^{(n)};t) dT \\
= \frac{1}{1-t} \Big(\displaystyle \sum_{\substack{i :\\  \lambda_{i} + \lambda_{2n} = 0}} t^{i-1} \Big) \int R_{\widehat{\lambda}, \widehat{w}}(x^{(n-1)}, y^{(n-1)};t) \Delta(x^{(n-1)};y^{(n-1)};t) dT  
\end{multline*}
where $\widehat{w}$ is $w$ with $y_{1}, x_{n}$ deleted and $\widehat{\lambda}$ is $\lambda$ with $\lambda_{i}$ and $\lambda_{2n}$ deleted (where the index $i$ is such that $\lambda_{i} + \lambda_{2n} = 0$).
\end{claim}
The proof is analogous to the proof of Claim \ref{conj3cl2} of the previous theorem.

Thus,
\begin{align*}
\int_{T} R_{\mu \bar{\nu}, w}(x^{(n)}, y^{(n)};t) \Delta(x^{(n)};y^{(n)};t) dT &= 0
\end{align*}
unless $\mu = \nu$.  Moreover, if $\mu = \nu$, the integral is
\begin{align*}
 \frac{1}{(1-t)^{l(\mu)}} v_{\mu+}(t) \int R_{0^{2n-2l(\mu)}, \delta}(x^{(n-l(\mu))}, y^{(n-l(\mu))};t)\Delta(x^{(n-l(\mu))}; y^{(n-l(\mu))};t) dT 
 \end{align*}
if $w = y_{1} \dots y_{l(\mu)} \delta x_{n-l(\nu)+1} \dots x_{n}$ for some permutation $\delta$ of $\{ y_{l(\mu)+1}, \dots, y_{n}, x_{1}, \dots, x_{n-l(\nu)} \}$ and $0$ otherwise.  

By Claim \ref{conj5norm}, we have
\begin{equation*}
\int_{T} R_{0^{2n-2l(\mu)}}(x^{(n-l(\mu))},y^{(n-l(\mu))};t) \frac{\Delta(x^{(n-l(\mu))};y^{(n-l(\mu))};t)}{\Big((2n-2l(\mu))!\Big)^{2}} dT = \frac{\phi_{2n-2l(\mu)}(t)}{(1-t)^{2n-2l(\mu)}\phi_{n-l(\mu)}(t^{2})}
\end{equation*}

Thus,
\begin{multline*}
\frac{1}{Z}\int_{T} P_{\mu \bar{\mu}}(x^{(n)},y^{(n)};t)\frac{1}{(n!)^{2}} \Delta(x^{(n)};y^{(n)};t) dT = \frac{\phi_{n}(t^{2}) }{v_{\mu+}(t)^{2}v_{(0^{2n-2l(\mu)})}(t)}\frac{v_{\mu+}(t)}{(1-t)^{l(\mu)}} \frac{\phi_{2n-2l(\mu)}(t)}{(1-t)^{2n-2l(\mu)}\phi_{n-l(\mu)}(t^{2})} \\
= \frac{(1-(t^{2})^{n-l(\mu)+1}) \cdots (1-(t^{2})^{n})}{v_{\mu+}(t)(1-t)^{2n-l(\mu)}}\frac{\phi_{2n-2l(\mu)}(t)}{v_{(0^{2n-2l(\mu)})}(t)} 
=\frac{(1-(t^{2})^{n-l(\mu)+1}) \cdots (1-(t^{2})^{n})}{v_{\mu+}(t) (1-t)^{l(\mu)}}.
\end{multline*}
where the last equality follows from the definition of $v_{(0^{2n-2l(\mu)})}(t)$.
Finally, one can check from the definition of the $C$-symbols that
\begin{align*}
C^{+}_{\mu}(t^{2n-2}t;0,t) &= 1 \\
C^{0}_{\mu}(t^{n}, -t^{n};0,t) &= \prod_{1 \leq i \leq l(\mu)} (1-t^{2(n+1-i)}) \\
C^{-}_{\mu}(t;0,t) &= (1-t)^{l(\mu)}v_{\mu+}(t).
\end{align*}
so that our formula gives
\begin{align*}
\frac{C^{0}_{\mu}(t^{n},-t^{n};0,t)}{C^{-}_{\mu}(t;0,t) C^{+}_{\mu}(t^{2n-2}t;0,t)},
\end{align*}
as desired.
\end{proof}

\begin{theorem}
(see \cite[Theorem 4.4]{RV}) Let $\lambda$ be a weight of the double cover of $GL_{2n}$, i.e., a half-integer vector such that $\lambda_{i} - \lambda_{j} \in \mathbb{Z}$ for all $i,j$.  Then 
\begin{align*}
\frac{1}{Z}\int P_{\lambda}^{(2n)}( \cdots t^{\pm 1/2}z_{i} \cdots ;t) \frac{1}{n!} \prod_{1 \leq i<j \leq n} \frac{(1-z_{i}/z_{j})(1-z_{j}/z_{i})}{(1-t^{2}z_{i}/z_{j})(1-t^{2}z_{j}/z_{i})} dT &=0,
\end{align*}
unless $\lambda = \mu \bar{\mu}$.  In this case, the nonzero value is
\begin{align*}
\frac{\phi_{n}(t^{2})}{(1-t)^{n}v_{\mu}(t) (1+t)(1+t^{2}) \cdots (1+t^{n-l(\mu)})}&= \frac{C^{0}_{\mu}(t^{n}, -t^{n};0,t)}{C^{-}_{\mu}(t;0,t) C^{+}_{\mu}(t^{2n-2}t;0,t)}.
\end{align*}  
\end{theorem}

\begin{proof}

As usual, note that $P_{\lambda}^{(2n)}( \cdots t^{\pm 1/2}z_{i} \cdots ;t)$ is a sum of $(2n)!$ terms, one for each permutation in $S_{2n}$.  We first note that many of these have vanishing integrals:

\begin{claim}
Let $w \in S_{2n}$ be a permutation of $(t^{\pm 1/2}z_{1}, \dots, t^{\pm 1/2}z_{n})$, such that for some $1 \leq i \leq n$ $\sqrt{t}z_{i}$ appears to the left of $\frac{z_{i}}{\sqrt{t}}$ in $w$.  Then 
\begin{align*}
\int R_{\lambda, w}^{(2n)}(\cdots t^{\pm 1/2}z_{i} \cdots ;t) \tilde \Delta_{S}^{(n)}(z;t^{2}) dT = 0.
\end{align*}
\end{claim}
To prove the claim note that $R_{\lambda, w}^{(2n)}( \cdots t^{\pm 1/2}z_{i} \cdots;t ) = 0$ in this case.  Indeed, we have the term
\begin{align*}
\frac{\sqrt{t}z_{i} - tz_{i}/\sqrt{t}}{\sqrt{t}z_{i} - z_{i}/\sqrt{t}} = \frac{tz_{i} - tz_{i}}{z_{i}(t-1)} =0
\end{align*}
appearing in the product defining the Hall--Littlewood polynomial.

Thus, we may restrict our attention to those permutations $w$ with $z_{i}/\sqrt{t}$ to the left of $\sqrt{t}z_{i}$ for all $1 \leq i \leq n$.  Moreover, we may order the variables so that $z_{i}/\sqrt{t}$ appears to the left of $z_{j}/\sqrt{t}$ for all $1 \leq i<j \leq n$. 
We compute the normalization first.
\begin{claim}\label{8.3.2norm}
We have 
\begin{align*}
Z = \int_{T} P_{0^{2n}}^{(2n)}(\cdots t^{\pm 1/2}z_{i} \cdots ;t) \frac{1}{n!}\tilde\Delta_{S}^{(n)}(z;t^{2}) dT = \frac{1}{v_{(0^{n})}(t^{2})} = \frac{(1-t^{2})^{n}}{(1-t^{2})(1-t^{4}) \cdots (1-t^{2n})}.
\end{align*}
\end{claim}
The proof follows by noting that $P_{0^{2n}}^{(2n)}( \cdots t^{\pm 1/2}z_{i} \cdots ;t) = 1$ and applying Theorem \ref{orthog}.

\begin{claim}
Let $w \in S_{2n}$ be a permutation with $z_{i}/\sqrt{t}$ to the left of $\sqrt{t}z_{i}$ for all $1 \leq i \leq n$ and $z_{i}/\sqrt{t}$ to the left of $z_{j}/\sqrt{t}$ for all $1 \leq i<j \leq n$, and $\sqrt{t}z_{1}$ in position $k$ for some $2 \leq k \leq 2n$.  Then
\begin{multline*}
\int_{T} R_{\lambda, w}^{(2n)}(\cdots t^{\pm 1/2}z_{i} \cdots ;t) \tilde\Delta_{S}^{(n)}(z;t^{2}) dT \\
= 
 \chi_{\lambda_{1} + \lambda_{k}=0} (1+t) t^{2n-k} \int_{T} R_{\widehat{\lambda}, \widehat{w}}^{(2(n-1))}( \cdots t^{\pm 1/2}z_{i} \cdots ;t) \tilde\Delta_{S}^{(n-1)}(z;t^{2}) dT  
\end{multline*}
where $\widehat{w}$ is the permutation $w$ with $z_{1}/\sqrt{t}$ and $\sqrt{t}z_{1}$ deleted, and $\widehat{\lambda}$ is the partition $\lambda$ with parts $\lambda_{1}$ and $\lambda_{k}$ deleted.
\end{claim}
To prove the claim, integrate with respect to $z_{1}$.  Note that if $\lambda_{1} + \lambda_{k} > 0$, the integral vanishes.  If $\lambda_{1} + \lambda_{k} < 0$, note that $\lambda_{2n} + \lambda_{j} < 0$ for all $1 \leq j \leq 2n-1$.  Integrate with respect to the last variable in $w$, and invert all variables to find the integral vanishes, as desired.

The above claim implies that the integral $\int_{T} R_{\lambda, w}^{(2n)}( \cdots t^{\pm 1/2}z_{i} \cdots ;t) \tilde\Delta_{S}^{(n)}(z;t^{2}) dT$ vanishes unless $\lambda = \mu \bar{\mu}$ for some $\mu$.  Moreover, if $\lambda = \mu \bar{\mu}$, the term integral vanishes unless
\begin{align*}
w( \cdots t^{\pm 1/2}z_{i} \cdots )^{\lambda}
\end{align*}
is a constant in $t$ (i.e., independent of $z_{i}$).
Thus, in the case $\lambda = \mu \bar{\mu}$, a computation gives that the total integral
\begin{multline*}
\int_{T} R_{\lambda}^{(2n)}( \cdots t^{\pm 1/2}z_{i} \cdots ;t)\frac{1}{n!} \tilde\Delta_{S}^{(n)}(z;t^{2}) dT \\
= (1+t)^{l(\mu)}v_{\mu+}(t) \int_{T} R_{0^{2(n-l(\mu))}}^{(2(n-l(\mu)))}(\cdots t^{\pm 1/2}z_{i} \cdots;t) \frac{1}{(n-l(\mu))!} \tilde\Delta_{S}^{(n-l(\mu))}(z;t^{2}) dT \\
= (1+t)^{l(\mu)}v_{\mu+}(t) \frac{(1-t^{2})^{n-l(\mu)}}{(1-t^{2})(1-t^{4}) \cdots (1-t^{2(n-l(\mu))})}v_{(0^{2(n-l(\mu))})}(t).
\end{multline*}
Multiplying this by $1/Zv_{\lambda}(t) = 1/Zv_{\mu+}(t)^{2}v_{(0^{2(n-l(\mu))})}(t)$ and simplifying gives the result.
\end{proof}

\begin{theorem}
(see \cite[Corollary 4.7(ii)]{RV}) Let $\lambda$ be a partition with $l(\lambda) \leq n$.  Then the integral
\begin{align*}
\int P_{\lambda}(x_{1}, \dots, x_{n};t^{2}) P_{m^{n}}(x_{1}^{-1}, \dots, x_{n}^{-1};t) \frac{1}{n!} \tilde \Delta_{S}^{(n)}(x;t) dT
\end{align*}
vanishes unless $\lambda = (2m)^{n} - \lambda$.
\end{theorem}
Note that the above integral gives the coefficient of $P_{m^{n}}(x;t)$ in the expansion of $P_{\lambda}(x;t^{2})$ as Hall--Littlewood polynomials with parameter $t$.  

\begin{proof}
Since $P_{m^{n}}(x_{1}^{-1}, \dots, x_{n}^{-1};t) = (x_{1}^{-1} \cdots x_{n}^{-1})^{m}$, an equivalent statement is the following:

Let $\lambda$ be a weight of $GL_{n}$ with possibly negative parts.  Then the integral
\begin{align*}
\frac{1}{Z} \int P_{\lambda}(x_{1}, \dots, x_{n};t^{2}) \frac{1}{n!} \tilde \Delta_{S}^{(n)}(x;t) dT
\end{align*}
vanishes unless $\lambda = \mu \bar{\mu}$, and in this case it is
\begin{align*}
\frac{(1-t^{n-2l(\mu)+1}) \cdots (1-t^{n}) t^{|\mu|}}{(1-t^{2})^{l(\mu)} v_{\mu+}(t^{2})}.
\end{align*}

We first compute the normalization $Z = \frac{1}{n!}\int P_{0^{n}}^{(n)}(x;t^{2}) \tilde \Delta_{S}^{(n)}(x;t) dT$.  Note that $P_{0^{n}}(x;t^{2}) = 1$, so we have
\begin{multline*}
Z = \frac{1}{n!}\int \tilde \Delta_{S}^{(n)}(x;t) dT 
= \frac{1}{n!} \int P_{0^{n}}^{(n)}(x;t) P_{0^{n}}^{(n)}(x^{-1};t) \tilde \Delta_{S}^{(n)} dT 
= \frac{1}{n!} \frac{n!}{v_{(0^{n})}(t)} 
= \frac{(1-t)^{n}}{(1-t)(1-t^{2}) \cdots (1-t^{n})}
\end{multline*}
using Theorem \ref{orthog}.  

Now we look at $\frac{1}{n!}\int R_{\lambda}(x_{1}, \dots, x_{n}; t^{2}) \tilde \Delta_{S}^{(n)}(x;t) dT$, which is a sum of $n!$ integrals---one for each $w \in S_{n}$.  By symmetry we have
\begin{align*}
\frac{1}{n!}\int R_{\lambda}^{(n)}(x; t^{2}) \tilde \Delta_{S}^{(n)}(x;t) dT &= \int R_{\lambda, \text{id}}^{(n)}(x;t^{2}) \tilde \Delta_{S}^{(n)}(x;t) dT,
\end{align*}
so we may restrict to the case $w = \text{id}$.  We assume $\lambda_{1} >0$: note that if $\lambda_{1} \leq 0$ we have $\lambda_{n} < 0$ (we are assuming $\lambda$ is not the zero partition) and we can invert all variables and make a change of variables to reduce to the case $\lambda_{1} > 0$.  Then the integral restricted to terms in $x_{1}$ is
\begin{multline*}
 \int_{T_{1}} x_{1}^{\lambda_{1}} \prod_{j > 1} \frac{x_{1} - t^{2}x_{j}}{x_{1} - x_{j}} \prod_{j > 1} \frac{(x_{1} - x_{j})(x_{j}-x_{1})}{(x_{1}-tx_{j})(x_{j}-tx_{1})} \frac{dx_{1}}{2\pi \sqrt{-1}x_{1}} 
= \int_{T_{1}} x_{1}^{\lambda_{1}} \prod_{j>1} \frac{(x_{1} - t^{2}x_{j})(x_{j}-x_{1})}{(x_{1}-tx_{j})(x_{j}-tx_{1})}\frac{dx_{1}}{2\pi \sqrt{-1}x_{1}} \\
= \sum_{j > 1} \frac{t^{\lambda_{1}}(1-t)^{2}}{(1-t^{2})} x_{j}^{\lambda_{1}} \prod_{i \neq 1,j} \frac{(tx_{j} - t^{2}x_{i})(x_{i}-tx_{j})}{(tx_{j}-tx_{i})(x_{i}-t^{2}x_{j})} 
= \sum_{j > 1} \frac{t^{\lambda_{1}}(1-t)^{2}}{(1-t^{2})}x_{j}^{\lambda_{1}} \prod_{i \neq 1,j} \frac{(x_{j} - tx_{i})(x_{i}-tx_{j})}{(x_{j}-x_{i})(x_{i}-t^{2}x_{j})}, 
\end{multline*}
where the second line follows by evaluating the residues at $x_{1} = tx_{j}$ for $j>1$.  
For each $j > 1$, we can combine this with the terms in $x_{j}$ from the original integrand.  The integral restricted to terms in $x_{j}$ is
\begin{multline*}
 \frac{t^{\lambda_{1}}(1-t)^{2}}{(1-t^{2})}\int_{T_{1}} x_{j}^{\lambda_{1}} \prod_{i \neq 1,j} \frac{(x_{j} - tx_{i})(x_{i}-tx_{j})}{(x_{j}-x_{i})(x_{i}-t^{2}x_{j})} x_{j}^{\lambda_{j}}\prod_{1 \neq i<j} \frac{x_{i} - t^{2}x_{j}}{x_{i}-x_{j}} \prod_{j<i} \frac{x_{j}-t^{2}x_{i}}{x_{j}-x_{i}}\\ 
\cdot \prod_{i \neq 1,j} \frac{(x_{i}-x_{j})(x_{j}-x_{i})}{(x_{i}-tx_{j})(x_{j}-tx_{i})} \frac{dx_{j}}{2\pi \sqrt{-1}x_{j}} 
= \frac{t^{\lambda_{1}}(1-t)^{2}}{(1-t^{2})} \int x_{j}^{\lambda_{1} + \lambda_{j}} (-1)^{n-j} \prod_{j < i} \frac{x_{j}-t^{2}x_{i}}{x_{i} - t^{2}x_{j}} \frac{dx_{j}}{2 \pi \sqrt{-1} x_{j}}.
\end{multline*}
Now, this is $0$ if $\lambda_{1} + \lambda_{j} > 0$ and 
\begin{align*}
\frac{t^{\lambda_{1}}(1-t)(t^{2})^{n-i}}{(1+t)}
\end{align*}
if $\lambda_{1} + \lambda_{j} = 0$.  Finally, if $\lambda_{1} + \lambda_{j} < 0$ note that $\lambda_{n} + \lambda_{i} < 0$ for all $1 \leq i < n$.  We can invert all variables and make a change of variables to arrive at the case $\lambda_{1} + \lambda_{j} > 0$, so the integral is zero by the above argument.  

Iterating this argument shows that the partition $\lambda$ must satisfy $\lambda_{i} + \lambda_{n+1-i} = 0$ for the integral to be nonvanishing.  Thus $\lambda = \mu \bar{\mu}$ for some $\mu$.  In this case, we compute from the above remarks:
\begin{multline*}
 \frac{1}{Z} \int P_{\lambda}^{(n)}(x;t^{2}) \frac{1}{n!} \tilde \Delta_{S}^{(n)}(x;t) dT 
=   \frac{1}{Z} \frac{1}{v_{\lambda}(t^{2})} \int R_{\lambda, \text{id}}^{(n)}(x;t^{2}) \tilde \Delta_{S}^{(n)}(x;t) dT \\
= \frac{\phi_{n}(t)}{(1-t)^{n}} \frac{t^{|\mu|}}{v_{\mu+}(t^{2})^{2} v_{(0^{n-2l(\mu)})}(t^{2})}  \frac{(1-t)^{l(\mu)}}{(1+t)^{l(\mu)}} v_{\mu+}(t^{2}) \int R_{0^{n-2l(\mu)}}(x;t^{2}) \frac{1}{(n-2l(\mu))!}\tilde \Delta_{S}^{(n)}(x;t)dT.
\end{multline*}
Using the computation of $Z$, this is equal to
\begin{multline*}
 \frac{\phi_{n}(t)}{(1-t)^{n}} \frac{t^{|\mu|}}{v_{\mu+}(t^{2})}  \frac{(1-t)^{l(\mu)}}{(1+t)^{l(\mu)}} \int P_{0^{n-2l(\mu)}}^{(n-2l(\mu))}(x;t^{2}) \frac{1}{(n-2l(\mu))!}\tilde \Delta_{S}^{(n)}(x;t) dT\\
= \frac{\phi_{n}(t)}{(1-t)^{n}} \frac{t^{|\mu|} }{v_{\mu+}(t^{2})} \frac{(1-t)^{l(\mu)}}{(1+t)^{l(\mu)}} \frac{(1-t)^{n-2l(\mu)}}{\phi_{n-2l(\mu)}(t)} 
= \frac{\phi_{n}(t)}{\phi_{n-2l(\mu)}(t)} \frac{t^{|\mu|}}{(1-t^{2})^{l(\mu)}v_{\mu+}(t^{2})} \\
= \frac{(1-t^{n-2l(\mu)+1}) \cdots (1-t^{n}) t^{|\mu|}}{(1-t^{2})^{l(\mu)}v_{\mu+}(t^{2})},
\end{multline*}
as desired.
\end{proof}

\end{document}